\documentclass{article}

\usepackage{microtype}
\usepackage{graphicx}
\usepackage{booktabs} % for professional tables
\usepackage{subcaption}

\usepackage{hyperref}

\usepackage{amsmath}        %alignment in equation
\usepackage{amsthm}
\usepackage{amsfonts}       % blackboard math symbols
\usepackage{amssymb}

\newcommand{\argmin}{\arg\!\min} 

\newcommand{\vect}{\mathbf}

\usepackage{mathtools}
\usepackage{commath}
\let\norm\undefined % <-- "Undefine" \norm
\DeclarePairedDelimiter\norm{\lVert}{\rVert}

\newtheorem{definition}{Definition}
\newtheorem{theorem}{Theorem}

\newtheorem{lemma}[theorem]{Lemma}
\newtheorem{assumption}{Assumption}
\usepackage[dvipsnames]{xcolor}

\newtheorem{innercustomgeneric}{\customgenericname}
\providecommand{\customgenericname}{}
\newcommand{\newcustomtheorem}[2]{%
  \newenvironment{#1}[1]
  {%
   \renewcommand\customgenericname{#2}%
   \renewcommand\theinnercustomgeneric{##1}%
   \innercustomgeneric
  }
  {\endinnercustomgeneric}
}

\DeclarePairedDelimiter{\ceil}{\lceil}{\rceil}

\newcustomtheorem{customthm}{Theorem}
\newcustomtheorem{customlemma}{Lemma}

\usepackage[shortlabels,inline]{enumitem}

\usepackage[accepted]{icml2019}

\icmltitlerunning{Combining Stochastic Adaptive Cubic Regularization and Negative Curvature for Nonconvex Optimization}

\begin{document}

\twocolumn[
% \icmltitle{Stochastic Adaptive Cubic Regularization with Exploiting Negative Curvature for Nonconvex Optimization}
\icmltitle{Combining Stochastic Adaptive Cubic Regularization and Negative Curvature for Nonconvex Optimization}



\icmlsetsymbol{equal}{*}

\begin{icmlauthorlist}
\icmlauthor{Seonho Park}{UFISE}
\icmlauthor{Seung Hyun Jung}{KITECH}
\icmlauthor{Panos M. Pardalos}{UFISE}
\end{icmlauthorlist}

\icmlaffiliation{UFISE}{Center for Applied Optimization, Department of Industrial and Systems Engineering, University of Florida, Gainesville, Florida, USA}
\icmlaffiliation{KITECH}{Construction Equipment Technology Center, Korea Institute of Industrial Technology (KITECH), Republic of Korea}

\icmlcorrespondingauthor{Seonho Park}{seonhopark@ufl.edu}

\icmlkeywords{Machine Learning, ICML, Stochastic Optimization, Adaptive Cubic Regularization, Negative Curvature}

\vskip 0.3in
]




\begin{abstract}
We focus on minimizing nonconvex finite-sum functions that typically arise in machine learning problems. 
In an attempt to solve this problem, the adaptive cubic regularized Newton method has shown its strong global convergence guarantees and ability to escape from strict saddle points. 
This method uses a trust region-like scheme to determine if an iteration is successful or not, and updates only when it is successful. 

In this paper, we suggest an algorithm combining negative curvature with the adaptive cubic regularized Newton method to update even at unsuccessful iterations. 
We call this new method Stochastic Adaptive cubic regularization with Negative Curvature (SANC). 
Unlike the previous method, in order to attain stochastic gradient and Hessian estimators, the SANC algorithm uses independent sets of data points of consistent size over all iterations.
It makes the SANC algorithm more practical to apply for solving large-scale machine learning problems. 
To the best of our knowledge, this is the first approach that combines the negative curvature method with the adaptive cubic regularized Newton method. 
Finally, we provide experimental results including neural networks problems supporting the efficiency of our method.
\end{abstract}

\section{Introduction}\label{sec:introduction}
We have focused on solving the nonconvex unconstrained problems.
Especially, in machine learning problems with large dataset, we frequently exploit stochastic optimization schemes; when it comes to solving an empirical risk minimization(ERM) problem, instead of calculating exact gradient and Hessian at an iterate, their estimators are used to find an update step. 
In this context, in order to find second order critical point, two of the most popular optimization approaches using Hessian are considered:
\begin{enumerate*}[label=(\roman*)]
\item the negative curvature method and
\item the Newton method.
\end{enumerate*}

Given an iterate, the negative curvature methods \cite{curtis2017exploiting,liu2018adaptive, cano2017using} need to calculate the eigenvector corresponding to the left-most eigenvalue of Hessian.
Using this eigenvector along with a novel approach to switch negative curvature direction and gradient descent direction, this method shows its convergence to a second order critical point.
% Depending on a novel measure of each negative curvature variant, this method alternates between negative curvature and gradient descent or combines both update steps.
We remark that because calculating exact eigenvector and corresponding left-most eigenvalue are time-consuming, this method may use approximate approaches such as the Lanczos method or the Oja's algorithm \cite{oja1982simplified}.

% Next Paragraph - Newton methods

The second approach is the Newton method. \cite{martens2010deep,martens2011learning,agarwal2017second,vinyals2012krylov}
Given an iterate, to find an update step, Newton method solves a linear system which denotes the first order necessary condition of optimality for a second-order Taylor approximation of an objective function.
To prevent the storage of the Hessian matrix which can be prohibitive to store especially for large-scale problems, We usually use Krylov subspace methods to iteratively attain an approximated solution of the linear system over a Krylov subspace.
The Krylov subspace methods also utilize Hessian-vector products \cite{pearlmutter1994fast} to construct a Krylov subspace whose dimension increases linearly over iterations.

In a similar vein, the cubic regularized Newton method \cite{nesterov2006cubic} has been used to escape from a strict saddle point and ultimately converge to a second-order stationary point for nonconvex problems. 
% Given an iterate $\vect{x}_t$, with a proper cubic coefficient $\sigma>0$, a \textit{local cubic model} $m_t(\vect{s})$ can be constructed as 
% \begin{equation}\label{eq_local_cubic_model1}
%     m_t(\vect{s}) := f(\vect{x}_t)+\nabla f(\vect{x}_t)^T\vect{s}+\frac{1}{2}\vect{s}^T\vect{H}(\vect{x}_t)\vect{s}+\frac{1}{3}\sigma\norm{\vect{s}}^3.
% \end{equation}
% The Newton step $\vect{s}_t$ can be obtained by minimizing Eq. \ref{eq_local_cubic_model1} over whole space or the Krylov subspace.
% % Also, approximate minimizer over the Krylov subspace of Eq. \ref{eq_local_cubic_model1} can be considered.
More recently, Adaptive Regularization algorithm using Cubics (ARC) and its stochastic variant \cite{cartis2011adaptive,cartis2011adaptive2,kohler2017sub,bergou2018line} are proposed. 
The ARC algorithm adopts the trust-region concept to determine if the obtained Newton update step is successful or not.
The next iterate is updated with the update step only when the update step is successful.

We focus our attention on integrating the negative curvature method into the ARC algorithm to update with a negative curvature direction when update Newton step is not successful .
% With this, the next iterate is updated with some direction even when the update step is not successful.
% Our method utilizes a direction of negative curvature to update an iterate even when an update Newton step $\vect{s}_t$ is not reliable.
Also, it is noted that calculating a negative curvature direction is efficient because it is on the same Krylov subspace where the update step of ARC is obtained.

% Contribution of this research
Thus, we make the following contributions in this paper.
\begin{itemize}
    \item We provide an optimization algorithm which mainly uses the ARC algorithm and also exploit negative curvature when the obtained update step is not successful.
    \item We also provide the theoretical results which include worst-case iteration complexity for achieving approximate first- and second-order optimality and corresponding lower bounds on the cardinalities of the data subsets on which stochastic gradient and Hessian estimators are calculated.
    \item We verify the efficiency of our method numerically by comparing with other standard first- and second-order methods for nonconvex optimization to show its prowess especially for neural networks problems.
\end{itemize}

\section{Related Works}\label{sec:related_works}
\textbf{Second-order methods}
In the stochastic optimization setting, second-order methods have been researched for decades.
J. Martens \cite{martens2010deep,martens2011learning} published papers about the Hessian-free optimization method to train neural networks including deep auto-encoder and recurrent neural networks.
% He argued in his papers that the use of a Gauss-Newton approximation instead of a Hessian approximation is more effective to train neural networks.
In these papers, it is argued that using a Gauss-Newton approximation instead of a Hessian approximation is more effective to train neural networks.
Also, he employs a damping coefficient to maintain the positive definiteness of a Hessian approximation, which helps restrict a step size.
Consequently, with a damping coefficient $\lambda$, a linear system is solved, $(\vect{H}+\lambda I)\vect{s}=-\vect{g}$, to attain an update step $\vect{s}$ by using the conjugate gradient algorithm, at every iterate.
% TODO: Also, a damping coefficient to maintain the positive definiteness of a Hessian approximation is employed, which helps to restrict a step size. Consequently, with a damping coefficient $\lambda$, a linear system is solved, $(\vect{H}+\lambda I)\vect{s}=-\vect{g}$, to attain an update step $\vect{s}$ by using the conjugate gradient algorithm, at every iterate.
The idea to use a damping coefficient is highly similar to that of adopting an adaptive cubic regularization term.
But it is noted that adaptive cubic regularization has a theoretical adaptive rule to guarantee its convergence whereas the Newton method with a damping coefficient does not provide.
% TODO: But it is noted that adaptive cubic regularization has a theoretical adaptive rule to guarantee its convergence whereas the Newton method with a damping coefficient does not. (마지막에 provide 만 뺀거임)

To circumvent solving a linear system, many researchers have considered the quasi-Newton method \cite{wang2017stochastic,byrd2016stochastic} as an alternative.
In stochastic optimization, they applied an L-BFGS update formula using stochastic gradient and Hessian estimators.
% TODO: In this paper, an L-BFGS update formula using stochastic gradient and Hessian estimators is applied for stochastic optimization.
However, since an assumption of positive definiteness of a Hessian is required to prove its convergence, it is not capable to provide convergence guarantees for nonconvex optimization problems.
Experience has shown that some gains in performance in machine learning applications can be achieved, but the full potential of the stochastic quasi-Newton schemes is not yet known. \cite{bottou2018optimization}

\textbf{Cubic regularized Newton method}
The cubic regularized Newton method is proposed by Griewank \cite{griewank1981modification} firstly and is also studied by Nesterov \textit{et al.} \cite{nesterov2006cubic}.
In the deterministic setting, the use of a cubic regularization as an over-estimator of an objective function to solve unconstrained optimization problems has been studied.
% They constructed and solved a second-order Taylor approximation with a cubic regularization term at every iterate.
A second-order Taylor approximation with a cubic regularization term is constructed and solved at each iteration.
% TODO: In this work, a second-order Taylor approximation with a cubic regularization term is constructed and solved at each iteration.
Nesterov and Polyak \cite{nesterov2006cubic} provided the theoretical analysis of its global convergence rate for nonconvex unconstrained problems.

Their theoretical results for the cubic regularization paved the way for the upcoming ARC algorithm \cite{cartis2011adaptive,cartis2011adaptive2}.
This algorithm uses an adaptive coefficient in the cubic regularization term.
This coefficient varies according to the difference between the local cubic model value and actual function value at the next iterate.
The way of adjusting a cubic coefficient in the ARC algorithm is similar to that of adjusting a radius of the trust region methods.
They also showed the mathematical proof for the global worst-case complexity bounds.
% TODO: Also, the mathematical proof for the global worst-case complexity bounds is given.

More recently, Kohler and Lucchi \cite{kohler2017sub} provided stochastic ARC variant for solving finite-sum structure nonconvex objectives.
They proved the lower bound on cardinalities of the subsampling subsets to calculate stochastic gradient and Hessian estimators.
% TODO: In this paper, the lower bound on cardinalities of the subsampling subsets to calculate stochastic gradient and Hessian estimators is proved.
But in their method, the sizes of the subsets are increasing over iterations which hinders training large scale neural networks.
% TODO: But in this method, the sizes of the subsets are increasing over iterations which hinders training large scale neural networks.

Wang \textit{et al.} \cite{wang2017stochastic,wang2018cubic} proposed approaches to incorporate a momentum acceleration, or a stochastic variance reduced gradient(SVRG)(\cite{johnson2013accelerating}) into a framework of the cubic regularized Newton method.
However, they did not use an adaptive cubic coefficient in their paper.
% TODO: However, an adaptived cubic coefficient did not used in this paper.

\textbf{Negative curvature approaches}
There have been relatively little works on the negative curvature methods.
Curtis and Robinson \cite{curtis2017exploiting} proposed several algorithms exploiting negative curvature for solving deterministic and stochastic optimization problems. 
In this work, a current iterate is updated with a direction of (stochastic) gradient descent and negative curvature.
The 'dynamic method' is also proposed which adaptively estimates a Lipschitz constant of gradient continuity to estimate a step length in the stochastic optimization framework.
% In stochastic optimization, they also proposed the 'dynamic method' which adaptively estimates a Lipschitz constant of gradient continuity to estimate a step length.
% TODO: In stochastic optimization, the 'dynamic method' is proposed which adaptively estimates a Lipschitz constant of gradient continuity to estimate a step length.
They showed that the 'dynamic method' is efficient in performance in stochastic optimization problems, but they did not provide a proof of its convergence.
% TODO: In this paper, it is shown that the 'dynamic method' is efficient in performance in stochastic optimization problems empirically, but a proof of its convergence is not provided.

Most recently, in a similar vein, Liu  \textit{et al.} \cite{liu2018adaptive} proposed the adaptive negative curvature descent (NCD) method which gives some adaptability to the terminating criteria of the Lanczos procedure depending on the magnitude of the subsampled gradients.
They also provided variants of the adaptive NCD method whose worst-case time complexity is $\mathcal{\widetilde{O}}(d/\epsilon^{3.5})$ for stochastic optimization, where $\mathcal{\widetilde{O}}$ hides logarithmic terms.

Carmon \textit{et al.} \cite{carmon2018accelerated} used negative curvature directions at the first phase of iterates and then switched it to accelerated stochastic gradient descent method when an iterate reaches an almost convex region.

Also, a direction of negative curvature has been used for escaping from a strict saddle point.
For the details on it, please refer to \cite{reddi2017generic,xu2018first}.

\section{Formulation of Stochastic Adaptive Cubic Regularization with Negative Curvature}\label{sec:SANC}
\subsection{Problem statement}
We consider the following stochastic optimization problem:
\begin{equation}\label{eq_stochastic_optimization}
    \min_{\vect{x} \in \mathbb{R}^d} f(\vect{x}) = \mathbb{E}\left[F(\vect{x},\xi) \right]
\end{equation}
where $F$ is a twice continuously differentiable function and possibly nonconvex, $\vect{x}\in \mathbb{R}^d$ is a variable, and $\mathbb{E}$ denotes the expectation with respect to $\xi$, a random variable with the distribution $P$.
If it is assumed that there is no prior knowledge about the distribution $P$, one may use the following finite-sum structure, as a proxy of the problem, of the form \ref{eq_stochastic_optimization}.
% TODO: We assume that we do not have any prior knowledge about the distribution $P$, thus, the following finite-sum structure of the form may be used as a proxy of the problem \ref{eq_stochastic_optimization}.
\begin{equation}\label{eq_finitesum}
    \min_{\vect{x}\in \mathbb{R}^d} f(\vect{x}) := \frac{1}{n}\sum_{i=1}^{n}f(\vect{x};\vect{q}_i,r_i)=\frac{1}{n}\sum_{i=1}^{n}f_i(\vect{x})
\end{equation}
where a function $f_i$ is an abbreviation of $f(\vect{x};\vect{q}_i,r_i)$, and $\{(\vect{q}_1,r_1),\dots,(\vect{q}_n,r_n)\}$ are $n$ sampled dataset.
This optimization problem \ref{eq_finitesum} is usually referred to as an empirical risk minimization (ERM) problem.
We frequently encounter this kind of optimization problems in supervised machine learning.
For example, given a dataset, where $\vect{q}_i \in \mathbb{R}^m$ is a feature and $r_i \in \mathbb{R}$ is its corresponding label for all $i \in \{1,\dots,n\}$, $f_i(\vect{x}) = 1-\tanh(r_i\vect{x}^T\vect{q}_i) + \lambda \norm{\vect{x}}_2^2$ is a nonconvex support vector machine (SVM) \cite{wangpardalos2017stochastic} with a convex regularization term.

Also, when it comes to online or stochastic optimization, it is assumed that it is not possible to access directly to an exact function value, gradient, or Hessian at $\vect{x}$, i.e., $f(\vect{x})$, $\nabla f(\vect{x})$, or $\vect{H}(\vect{x})$.
Instead, we may resort to using unbiased estimators of the gradient and Hessian
%TODO: Instead, we may resort to using stochastic gradient and Hessian estimator.
As a prototypical method to attain stochastic gradient and Hessian estimators, one frequently consider using the mini-batch approach.
We construct sets of data points independently drawn from the dataset, $\mathcal{S}_\vect{g}$, $\mathcal{S}_\vect{B}$.
The stochastic gradient and Hessian estimators, say, $\vect{g}(\vect{x})$ and $\vect{B}(\vect{x})$, can be defined as,
\begin{equation}\label{eq_unbiased_estimator}
\begin{aligned}
\vect{g}(\vect{x}) := \frac{1}{|\mathcal{S}_\vect{g}|}\sum_{i\in \mathcal{S}_\vect{g}}\nabla f_i(\vect{x})\\
\vect{B}(\vect{x}) := \frac{1}{|\mathcal{S}_\vect{B}|}\sum_{i\in \mathcal{S}_\vect{B}}\nabla^2 f_i(\vect{x})
\end{aligned}
\end{equation}
In this context, we aim to find an $\epsilon$-approximate second-order stationary point $\vect{x}^*$ of the problem \ref{eq_finitesum} which satisfies,
\begin{equation}\label{eq_second_order_stationary_point}
    \norm{\nabla f(\vect{x}^*)} \leq \epsilon \text{ and } \lambda_{\min}\left(H(\vect{x}^*) \right) \geq -\epsilon
\end{equation}
where $\lambda_{\min}(\cdot)$ denotes the left-most eigenvalue.

For brevity, we use $\norm{\cdot}$ as an abbreviation of the $l2$ norm, $\norm{\cdot}_2$, in what follows. For the same reason, we also use $\vect{g}_t$, and $\vect{B}_t$ instead of $\vect{g}(\vect{x}_t)$, and $\vect{B}(\vect{x}_t)$, respectively.

\subsection{Algorithm}\label{subsec:algorithm}

\begin{algorithm*}[ht]
   \caption{Stochastic Adaptive cubic regularization method with Negative Curvature (SANC)}   \label{alg:SANC}
\begin{algorithmic}[1]
   \STATE {\bfseries Input:} initial point $\vect{x}_0\in\mathbb{R}^d$, $\gamma>1$, $1>\eta_2>\eta_1>0$, $L_1>0$, $L_2>0$, and $\sigma_0>0$
   \FOR{$t=0,1,\dots$ {until convergence} } 
    \STATE Set $\mathcal{S}_{\vect{g},t}$ and $\mathcal{S}_{\vect{B},t}$ of the dataset
    \STATE Sample $\vect{g}_t := \frac{1}{|\mathcal{S}_{\vect{g},t}|}\sum_{i\in \mathcal{S}_{\vect{g},t}}\nabla f_i(\vect{x}_t)$
    \STATE Sample $\vect{B}_t := \frac{1}{|\mathcal{S}_{\vect{B},t}|}\sum_{i\in \mathcal{S}_{\vect{B},t}}\nabla^2 f_i(\vect{x}_t)$
    \STATE Compute $\vect{v}_t$ such that Eq. \ref{eq:v_condition} holds with $\epsilon'=\max\{\epsilon,\norm{\vect{g}_t}\}/2$
    % \STATE Compute 
    % \begin{equation}\label{eq:d_rule}
    %     \vect{d}_t=
    %     \begin{cases}
    %         -\frac{2|\vect{v}_t^T\vect{B}_t\vect{v}_t|}{L_2}z\vect{v}_t & \text{if } \frac{2(-\vect{v}_t^T\vect{B}_t\vect{v}_t)^3}{3L_2^2}-\frac{\epsilon (\vect{v}_t^T\vect{B}_t\vect{v}_t)^2}{6L_2^2} > \frac{\norm{\vect{g}_t}^2}{4L_1}-\frac{\epsilon_{\vect{g}}^2}{L_1}\\
    %         -\frac{1}{L_1}\vect{g}_t & \text{otherwise}
    %     \end{cases}
    % \end{equation}
    % where $z\in \{-1,1\}$ with equal probability is a so-called Rademacher random variable.
    \STATE Compute $\vect{s}_t$ by attaining $\vect{u}_t$ of Eq. \ref{eq:krylov_minimizer} and retrieving by Eq. \ref{eq:retrieving_s} such that Eq. \ref{eq:s_condition1} and Eq. \ref{eq:s_condition2} hold
    \STATE Compute $f(\vect{x}_t+\vect{s}_t)$ and 
        \begin{equation}\label{eq:rho}
            \rho_t = \frac{f(\vect{x}_t)-f(\vect{x}_t+\vect{s}_t)}{f(\vect{x}_t)-\widetilde{m}_t(\vect{s}_t)}
        \end{equation}
    \IF{$\rho_t \geq \eta_1$ [(very) successful iteration]}
        \STATE $\vect{x}_{t+1}=\vect{x}_t+\vect{s}_t $
    \ELSE 
        \STATE Compute $\vect{d}_t$ with Eq. \ref{eq:d_rule}
        \STATE $\vect{x}_{t+1}=\vect{x}_t+\vect{d}_t $
    \ENDIF
    % \STATE Update
    % \begin{equation}\label{eq:update_rule}
    %     \vect{x}_{t+1}=
    %     \begin{cases}
    %       \vect{x}_t+\vect{s}_t & \text{if } \rho_t \geq \eta_1 \text{ [(very) successful iteration]}\\    
    %       \vect{x}_t+\vect{d}_t & \text{otherwise }\text{ [unsuccessful iteration]} 
    %     \end{cases}
    % \end{equation} 
    \STATE Update
    \begin{equation*}\label{eq:rho_update_rule}
        \sigma_{t+1}=
        \begin{cases}
            \max\{\min\{\sigma_t,\norm{\vect{g}_t}\},\epsilon_m\} & \text{if } \rho_t>\eta_2  \text{ [very successful iteration]} \\
            \sigma_t & \text{if } \eta_2\geq \rho_t \geq \eta_1 \text{ [successful iteration]} \\
            \gamma\sigma_t & \text{otherwise } \text{ [unsuccessful iteration]} 
        \end{cases}
    \end{equation*}
    where $\epsilon_m$ is the machine precision.
   \ENDFOR
\end{algorithmic}
\end{algorithm*}

In this section, we describe a general approach to our method which is presented as Algorithm \ref{alg:SANC}.
We call this method Stochastic Adaptive cubic regularization with Negative Curvature (SANC) in what follows.
The main framework of Algorithm \ref{alg:SANC} is same as that of the ARC algorithm.
Given an iterate $\vect{x}_t$, the algorithm determines if the Newton step $\vect{s}_t$ is successful or not. 
It depends on the prescribed parameter $\eta_1$ and $\eta_2$.
It also depend on a measure $\rho_t$ in Algorithm \ref{alg:SANC} which calculates the fraction of the predicted model value decrease $f(\vect{x}_t)-\widetilde{m}_t(\vect{s}_t)$ and the actual decrease $f(\vect{x}_t)-f(\vect{x}_t+\vect{s}_t)$.
If the difference between the predicted and the actual decrease is small enough, then the iteration may become \textit{successful}, otherwise, it is \textit{unsuccessful}.

If the iteration is successful, the next iterate is set to $\vect{x}_{t+1}=\vect{x}_t+\vect{s}_t$.
If the iteration is unsuccessful, we try to find a negative curvature direction to update. 

After updating the iterate, we need to update the cubic coefficient.
When the iteration is very successful, an adaptive cubic coefficient $\sigma_t$ is decreasing according to Eq. \ref{eq:rho_update_rule} in Algorithm \ref{alg:SANC}.
Also, if the iteration is unsuccessful, $\sigma_t$ is increasing, which leads to the Newton step with a small magnitude.
Thus, this adaptive cubic coefficient $\sigma_t>0$ plays a similar role as a reciprocal of radius in the trust region methods.

% As said before, we incorporate the negative curvature method into the ARC algorithm.
% Even at an unsuccessful iteration, a current iteration may have an opportunity to be updated with a negative curvature direction if a current iteration point lies on a locally nonconvex region. 
We describe the details on how we attain a Newton step and a negative curvature direction.

\subsection{Finding Newton step from local cubic model}\label{subsec:find_newton_step}
In this subsection, we elucidate a way to calculate a Newton step $\vect{s}_t$.
Given an iterate $\vect{x}_t$, we can attain the stochastic gradient and Hessian estimator, $\vect{g}_t$, $\vect{B}_t$, over $\mathcal{S}_{\vect{g},t}$ and $\mathcal{S}_{\vect{B},t}$, respectively.
With $\vect{g}_t$, $\vect{B}_t$, and a sufficiently large $\sigma_t>0$, we establish an approximate local cubic model $\widetilde{m}_t(\vect{s})$ of the form,
\begin{equation}\label{eq:local_cubic_model2}
    \widetilde{m}_t(\vect{s}) := f(\vect{x}_t)+\vect{g}_t^T\vect{s}+\frac{1}{2}\vect{s}^T\vect{B}_t\vect{s}+\frac{1}{3}\sigma_t\norm{\vect{s}}^3
\end{equation}
We can attain a global minimizer $\vect{s}_t$ corresponding to Eq. \ref{eq:local_cubic_model2} by $\vect{s}_t := \argmin_{\vect{s}\in \mathbb{R}^d}{\widetilde{m}_t(\vect{s})}$ regardless of the positive definiteness of $\vect{B}_t$.
As Cartis \textit{et al.} \cite{cartis2011adaptive} pointed out, attaining the global minimizer of Eq. \ref{eq:local_cubic_model2} is time-consuming.
So they suggested that solving Eq. \ref{eq:local_cubic_model2} is relaxed by satisfying the Cauchy condition,
\begin{equation}\label{eq:cauchy_condition}
\widetilde{m}_t(\vect{s}_t)\leq \widetilde{m}_t(\vect{s}_t^c),\:\: \vect{s}_t^c=-\alpha_c\vect{g}_t
\end{equation}
where $\alpha_c \in \argmin_{\alpha>0}\widetilde{m}_t(-\alpha \vect{g}_t)=f(\vect{x}_t)-\alpha\norm{\vect{g}_t}^2+\frac{1}{2}\alpha^2\vect{g}_t^T\vect{B}_t\vect{g}_t+\frac{1}{3}\sigma_t\alpha^3\norm{\vect{g}_t}^3 $.
Also, they provided a convergence guarantee with the satisfaction of the Cauhy condition.

As a consequence, the use of the Lanczos method can be considered to seek an approximate minimizer of Eq. \ref{eq:local_cubic_model2} over the Krylov subspace.
The Lanczos method is a kind of iterative methods to solve a linear system.
The Lanczos method constructs an orthogonal matrix, $\vect{Q}_t$, sequentially whose columns are a basis of a Krylov subspace $\mathcal{K}_j$ of $\mathbb{R}^d$.
It is noted that if the Krylov subspace is formed by $\vect{g}_t$ and $\vect{B}_t$ successively as
\begin{equation}\label{eq:krylov_subspace}
    \mathcal{K}_j(\vect{g}_t,\vect{B}_t)=span\left\{ \vect{g}_t,\vect{B}_t\vect{g}_t,\vect{B}_t^2\vect{g}_t,\dots,\vect{B}_t^{j-1}\vect{g}_t \right\},    
\end{equation}
then $\mathcal{K}_j(\vect{g}_t,\vect{B}_t)$ always contain a subspace $\vect{g}_t$ even at Lanczos iteration $j=1$.
Hence, an approximate minimizer $\vect{s}_t$ over the Krylov subspace $\mathcal{K}_{j\geq1}(\vect{g}_t,\vect{B}_t)$ always holds the Cauchy condition \ref{eq:cauchy_condition} when preconditioning is not applied.
As a consequence, the approximate minimizer $\vect{s}_t$ can be retrieved as 
\begin{equation}\label{eq:retrieving_s}
    \vect{s}_t = \vect{Q}_t\vect{u}_t   
\end{equation}
where 
\begin{equation}\label{eq:krylov_minimizer}
    \vect{u}_t = \argmin_{\vect{u}\in \mathcal{K}_j(\vect{g}_t,\vect{B}_t)}\gamma_0\vect{u}^Te_1+\frac{1}{2}\vect{u}^T\vect{T}\vect{u}+\frac{1}{3}\sigma_t\norm{\vect{u}}^3
\end{equation}
In Eq. \ref{eq:krylov_minimizer}, $\gamma_0=\norm{\vect{g}_t}$, $e_1$ is $1st$ column vector of identity matrix and $\vect{T}$ is a $j$ by $j$ symmetric tridiagonal matrix.
We would like to point out that solving Eq. \ref{eq:krylov_minimizer} is more efficient than solving Eq. \ref{eq:local_cubic_model2} because of its smaller dimension and a structured construction of $\vect{T}$.
Carmon and Duchi \cite{carmon2016gradient} argued that it takes $\mathcal{O}(\log(1/\epsilon))$ iterations with a gradient descent method to find an $\epsilon$-approximate second-order stationary point for Eq. \ref{eq:local_cubic_model2} with mild assumptions. 
Thus, it is obvious that solving Eq. \ref{eq:krylov_minimizer} requires much less iterations than solving Eq. \ref{eq:local_cubic_model2} because it uses a structured tridiagonal matrix and has smaller dimension.

As shown in Lemma 3.2 in \cite{cartis2011adaptive}, given an iterate $\vect{x}_t$, an global and approximate minimizer of $\widetilde{m}_t(\vect{s})$ over the Krylov subspace always satisfy the following requirements,
\begin{subequations}
\begin{align}
    \vect{g}_t^T\vect{s}_t+\vect{s}_t^T\vect{B}_t\vect{s}_t+\sigma_t\norm{\vect{s}_t}^3=0\label{eq:s_condition1}\\
    \vect{s}_t^T\vect{B}_t\vect{s}_t+\sigma_t\norm{\vect{s}_t}^3\geq0\label{eq:s_condition2}
\end{align}
\end{subequations}
These statements \ref{eq:s_condition1},\ref{eq:s_condition2} play a vital role to prove the proofs of the main theorem of our algorithm in Section \ref{sec:theoretical_analysis}.

It is not surprising that the Cauchy point $\vect{s}_t^c$ also satisfy Eq. \ref{eq:s_condition1} and Eq. \ref{eq:s_condition2}.
If we use a larger Krylov subspace, we can attain more accurate $\vect{s}_t$ which improves the performance of the algorithm. 
However, at the same time, at each Lanczos iteration, we need $\mathcal{O}\left(|\mathcal{S}_{\vect{B}}|d \right)$ operations so we have to choose the termination condition of the Lanczos iteration carefully.
In the next subsection, we describe details on how to attain a direction of negative curvature which can be utilized at unsuccessful iterations.

\subsection{Finding negative curvature}\label{subsec:negative_curvature_step}

Given an iterate $\vect{x}_t$, suppose that the stochastic Hessian estimator $\vect{B}_t$ is indefinite and $\vect{s}_t^T\vect{B}_t\vect{s}_t<0$. 
With a fixed $\sigma_t$, the magnitude of the calculated Newton step $\vect{s}_t$ should be relatively larger than that with an arbitrary positive definite Hessian because it is needed to satisfy Eq. \ref{eq:s_condition2}.
So, we can conjecture that if $\sigma_t$ is not large enough, then $\vect{x}_t+\vect{s}_t$ is not in the vicinity of $\vect{x}_t$.
As a consequence, it is likely that the Newton step $\vect{s}_t$ becomes unsuccessful.

At this point, if we know that the stochastic Hessian estimator $\vect{B}_t$ is indefinite, then we can utilize a negative curvature direction to update an iterate.
Luckily, from the previous subsection, the Newton step $\vect{s}_t$ is achieved by the Lanczos method, and if the orthogonal matrix $\vect{Q}_t$ is saved or regenerated after terminating the Lanczos iterations, we can utilize the following relationship,
\begin{equation}\label{eq:tridiagonal_relationship}
    \vect{Q}_t^T\vect{B}_t\vect{Q}_t=\vect{T}.
\end{equation}
It is noted that since it is structured, eigenpairs $(\lambda_i^{(j)},\vect{v}_i^{(j)})$ of $\vect{T}$ can be cheaply computed.
Also the left-most eigenvalue $\lambda_1^{(j)}$ and corresponding eigenvector $\vect{v}_1^{(j)}$ of $\vect{T}$ are used to approximate the left-most eigenpair of $\vect{B}_t$.
The approximate left-most eigenvector of $\vect{B}_t$ can be attained by, $\vect{v}_t=\vect{Q}_t\vect{v}_1^{(j)}$.
We can approximate the left-most eigenpair of $\vect{B}_t$ by that of $\vect{T}$, i.e., $ (\lambda_t,\vect{v}_t) = (\lambda_1^{(j)},\vect{Q}_t\vect{v}_1^{(j)})$.
This approximation is often called the Ritz approximation.

Please note that with a small $\epsilon'>0$, there exists $\vect{v}\in \mathbb{R}^d$ such that 
\begin{equation}\label{eq:v_condition}
 \lambda_{\min}\left( \vect{B}_t \right) \geq \vect{v}^T\vect{B}_t\vect{v}-\epsilon'   
\end{equation}
where $\vect{v}$ can be achieved within $\mathcal{\widetilde{O}}\left(\frac{d}{\sqrt{\epsilon'}}\right)$ by the Lanczos method, as proved in \cite{kuczynski1992estimating}.

Usually, the Lanczos iteration is truncated at $j\ll d$ especially for large-scale problems in order to make the algorithm efficient.
Also, computing eigenpairs of $\vect{T}$ takes only $\mathcal{O}(j^2)$ operations because $\vect{T}$ is structured. 
So, it is computationally cheap to calculate the left-most eigenpair of $\vect{T}$ in large-scale problems.
Also, one can use divide-and-conquer algorithm \cite{coakley2013fast} to calculate the eigenpair of $\vect{T}$.

In order to prove the convergence guarantee of our algorithm, we maintain a framework of the ARC algorithm and use a negative curvature step $\vect{d}_t$ at an unsuccessful iteration \cite{liu2018adaptive} as
\begin{equation}\label{eq:d_rule}
    \vect{d}_t=
    \begin{cases}
        -\frac{2|\vect{v}_t^T\vect{B}_t\vect{v}_t|}{L_2}z\vect{v}_t & \text{if } \frac{2(-\vect{v}_t^T\vect{B}_t\vect{v}_t)^3}{3L_2^2}-\frac{\epsilon (\vect{v}_t^T\vect{B}_t\vect{v}_t)^2}{6L_2^2} > \frac{\norm{\vect{g}_t}^2}{4L_1}-\frac{\epsilon_{\vect{g}}^2}{L_1}\\
        -\frac{1}{L_1}\vect{g}_t & \text{otherwise}
    \end{cases}
\end{equation}
where $z$ is the Rademacher random variable.
% seen in Eq. \ref{eq:d_rule} in Algorithm \ref{alg:SANC}.
In Eq. \ref{eq:d_rule}, $\frac{2(-\vect{v}_t^T\vect{B}_t\vect{v}_t)^3}{3L_2^2}-\frac{\epsilon (\vect{v}_t^T\vect{B}_t\vect{v}_t)^2}{6L_2^2}$ is the expected decrease when the direction is negative curvature, whereas $\frac{\norm{\vect{g}_t}^2}{4L_1}-\frac{\epsilon_{\vect{g}}^2}{L_1}$ is the expected decrease when stochastic gradient descent is used. 
In general, $\vect{d}_t$ is a combination of negative curvature and stochastic gradient descent which leads to decreasing objective a lot as possible.
Thus, the SANC updates at every iteration even if the Newton direction is not successful, which enables the algorithm to be more efficient than the ARC algorithm.
Now, we present the theorem regarding the convergence and the sampling approximations of gradients and Hessian in the next section.

\section{Theoretical Analysis}\label{sec:theoretical_analysis}
In this section, we provide the worst-case iteration complexity for achieving approximate first- and second-order optimality of the SANC algorithm. 
Firstly, we explain some assumptions regarding the gradient and Hessian approximations.
Also, we describe some assumptions regarding the convergence of the algorithm and the terminating criteria of the Lanczos method, which are crucial to the proof of the convergence of the SANC algorithm.
\subsection{Definition and assumptions}
\begin{definition}[Approximate Second-Order Stationary Point]\label{def:second_order_stationary_point}
With a small $\epsilon>0$, $x_{t^*}$ of iterates $\{x_t\}_{t\geq0}$ is called an $\epsilon$-approximate second-order stationary point if it satisfies
\begin{equation}
    \norm{\nabla f(\vect{x_{t^*}})} \leq \epsilon \text{ and } -\lambda_{\min}(\vect{Q}_{t^*}^T\vect{B}_{t^*}\vect{Q}_{t^*})\leq \epsilon.
\end{equation}
\end{definition}
Notably, even though the $\epsilon$-approximate second-order stationary point does not necessarily denote that it is close to an exact second-order stationary point, but if strict saddle property \cite{lee2016gradient} holds, then it is guaranteed that the $\epsilon$-approximate second-order stationary point is in the vicinity to the exact second-order stationary point for sufficiently small $\epsilon$.

To establish the convergence, we need to assume the following statements.
\begin{assumption}[Function Boundness]\label{assump:function_boundness}
Nonconvex function $f_i$ is twice differentiable and bounded below by $f_{inf}>-\infty$ for all $i$. 
Thus, nonconvex function $f$ is also twice differentiable and bounded below by $f_{inf}>-\infty$.
\end{assumption}
Assumption \ref{assump:function_boundness} plays an important role to impose an upper bound on the sum of the sufficient decreases of two consecutive function values.

\begin{assumption}[Lipschitz Continuity] \label{assump:Lipschitz}
For all $i$, the function $f_i$, $\nabla f_i$, and $\vect{H}_i$ are Lipschitz continuous with Lipschitz constants $L_0$, $L_1$, and $L_2$, respectively.
\end{assumption}
The assumption that $f_i$ is Lipschitz continuous is relatively rare to use, but in our analysis, it is used to bound the cardinality of the set for stochastic gradient estimator. 
Instead, we can also assume the exponential tail behavior of gradients realization.

To attain stochastic gradient and Hessian estimators, we resort to the following assumptions to independently draw the sets from the dataset.
\begin{assumption}[Gradient and Hessian Approximation Bounds]\label{assump:grad_hess_bound}
For all iteration $t\geq0$,
\begin{subequations}
\begin{align}
    \norm{\vect{g}_t - \nabla f(\vect{x}_t)} \leq \epsilon_{\vect{g}}\label{eq:grad_bound}\\
    \norm{\vect{B}_t - \vect{H}(\vect{x}_t)} \leq \epsilon_{\vect{B}}\label{eq:hessian_bound}
\end{align}
\end{subequations}
With sufficiently small $\epsilon_{\vect{g}}>0$ and $\epsilon_{\vect{B}}>0$
\end{assumption}

Assumption \ref{assump:grad_hess_bound} is different with those in \cite{cartis2011adaptive,kohler2017sub}. 
In their assumption, $\norm{\vect{s}_t}$ was included on the RHS of the inequality. 
As an iteration point is approaching an optimum, the discrepancy between exact and approximate gradient/Hessian needs a tighter agreement because $\norm{\vect{s}_t}$ tends to decrease.
Smaller $\norm{\vect{s}_t}$ makes the size of the sets increasing, which is harmful to apply for solving large-scale machine learning problems.

On the contrary, in our analysis, we can maintain the size of subsets consistent because of Assumption \ref{assump:grad_hess_bound}.
It makes our algorithm more effective to solve large-scale problems.

The following assumption about the termination criteria of the Lanczos method is used to establish the proof of Lemma \ref{lemma:lower_bound_s}.

\begin{assumption}[Lanczos Method Termination Criteria]\label{assump:lanczos_termination}
For each iteration $t$, let us assume that the Lanczos procedure stops at the Lanczos iteration $k$ when the following criterion holds,
\begin{equation}
    \norm{\nabla_s \widetilde{m}_t(\vect{s}_{t,k})}\leq \epsilon_\vect{g}
\end{equation}
\end{assumption}
The above assumption is used when we establish the proof of Lemma \ref{lemma:lower_bound_s}.
The proof uses the norm of discrepancy between $\nabla f(\vect{x}_t+\vect{s}_t)$ and $\nabla_s \widetilde{m}_t(\vect{s}_t)$.
Also, we need to assume the following statement to establish the upper bound on $\sigma_t$, which has a vital role in the convergence proof of our method.

\begin{assumption}[Convergence]\label{assump:convergence}
For Algorithm \ref{alg:SANC}, the process $\{\vect{x}_t\}_{t\geq0}$ stops when the following criteria hold, 
\begin{equation}
    \vect{d}_t = 0 \text{ and } \norm{\vect{s}_t} \leq \epsilon_s
\end{equation}
\end{assumption}

\subsection{Gradient and Hessian sampling bound}
Now, we establish the lower bound on the cardinalities of the sets for stochastic gradient and Hessian estimators.
The following theorems are based on the Bernstein inequalities in \cite{gross2011recovering}.
\begin{theorem}[Gradient Sampling Bound]\label{theorem:gradient_sampling_bound}
If
\begin{equation}\label{eq:grad_sampling_bound}
|\mathcal{S}_{\vect{g},t}| \geq  \frac{4L_0^2\left( 1+2\sqrt{\log{\frac{1}{\delta}}} \right)^2}{\epsilon_\vect{g}^2},
\end{equation}
then $\vect{g}_t$ satisfies the gradient approximation bound Eq. \ref{eq:grad_bound} for all $t\geq0$ with high probability $1-\delta$. 
\end{theorem}

\begin{theorem}[Hessian Sampling Bound]\label{theorem:hessian_sampling_bound}
If
\begin{equation}\label{eq:hessian_sampling_bound}
    |\mathcal{S}_{\vect{B},t}| \geq \frac{16L_1^2\log{\frac{2d}{\delta}}}{\epsilon_\vect{B}^2},
\end{equation}
then $\vect{B}_t$ satisfies the Hessian approximation bound Eq. \ref{eq:hessian_bound} for all $t\geq0$ with high probability $1-\delta$.
\end{theorem}
The proofs regarding Theorem \ref{theorem:gradient_sampling_bound}, \ref{theorem:hessian_sampling_bound} are in the supplementary material.
As shown in the Theorem \ref{theorem:gradient_sampling_bound}, \ref{theorem:hessian_sampling_bound}, the cardinalities of the sets $\mathcal{S}_\vect{g}$, $\mathcal{S}_\vect{B}$, do not depend on the norm of $\vect{s}$. 
Because these cardinalities bounds rely only on the constants, $L_1$, $d$, $\delta$, and $\epsilon_\vect{B}$, they can be consistent over all iterations.
Now, with the sampling bounds, we can establish the worst-case iteration complexity for achieving approximate first- and second-order optimality.

\subsection{Worst-case iteration complexity bound for achieving approximate first- and second-order optimality}

\begin{lemma}[Cubic Regularization Coefficient Bound]\label{lemma:cubicregularizationbound}
Let the A\ref{assump:Lipschitz}, A\ref{assump:grad_hess_bound} and A\ref{assump:convergence} hold. 
For finite positive values $\sigma_{\min}$ and $\sigma_{\max}$,
it holds
\begin{equation}\label{eq:sigma_limit}
    \sigma_{\min} \leq \sigma_t \leq \sigma_{\max} \: \text{ for all }t\geq0
\end{equation}
where $\sigma_{\max}:=\max\left\{\sigma_0,\gamma\left( \frac{3}{2}L_2+\frac{3(\epsilon_\vect{g}+\frac{1}{2}\epsilon_\vect{B})}{\epsilon_s^2} \right) \right\}$.
\end{lemma}

The following lemma is about the sufficient decrease with the negative curvature step $\vect{d}_t$.
\begin{lemma}[Sufficient Decrease with Negative Curvature Step, Lemma 3 in \cite{liu2018adaptive}]\label{lemma:sufficient_decrease_nc}
When $\vect{v}_t^T\vect{B}_t\vect{v}_t\leq -\epsilon'/2$ and $\epsilon_\vect{B} \leq \epsilon'/12$, the negative curvature step $\vect{d}_t$ satisfies that 
\begin{equation}\label{eq:d_update_bound}
    f(\vect{x}_t)-\mathbb{E}[f(\vect{x}_t+\vect{d}_t)] \geq \max\left\{ \frac{{\epsilon'}^3}{24L_2^2},\frac{\norm{\vect{g}_t}^2}{4L_1}-\frac{\epsilon_\vect{g}^2}{L_1} \right\}
\end{equation}
with high probability $1-\delta$.
\end{lemma}
In Lemma \ref{lemma:sufficient_decrease_nc}, we use the expectation of $f(\vect{x}_t+\vect{d}_t)$ with respect to the Rademacher random variable $z$.
This random variable takes the main role to simplify the lower bound of Eq. \ref{eq:d_update_bound}.
Before describing the sufficient decrease with the Newton step $\vect{s}_t$, the following lemmas are needed.

\begin{lemma}[Local Cubic Model Decrease, Lemma 3.3 in \cite{cartis2011adaptive}]\label{lemma:local_cubic_model_decrease}
Suppose that $\vect{s}_t$ satisfies Eq. \ref{eq:s_condition1} and Eq. \ref{eq:s_condition2}. Then for all (very) successful iterations $t$, 
\begin{equation}\label{eq:cubic_model_bound}
    f(\vect{x}_t)-\widetilde{m}_t(\vect{s}_t) \geq \frac{1}{6}\sigma_t\norm{\vect{s}_t}^3
\end{equation}
\end{lemma}

\begin{lemma}[Lower Bound of $\vect{s}_t$]\label{lemma:lower_bound_s}
Suppose that A\ref{assump:Lipschitz}, A\ref{assump:grad_hess_bound}, A\ref{assump:lanczos_termination}, and A\ref{assump:convergence} hold. 
For all (very) successful iterations $t$, the Newton step $\vect{s}_t$ satisfies that
\begin{equation}\label{eq:s_lower_bound}
    \norm{\vect{s}_t} \geq \frac{1}{\kappa_l}    \left(-\epsilon_\vect{B}+\sqrt{\epsilon_\vect{B}^2+2\kappa_l\left( \norm{\nabla f(\vect{x}_t+\vect{s}_t)}-2\epsilon_\vect{g} \right)} \right)
\end{equation}
where 
\begin{equation}
    \kappa_l:=L_2+2\sigma_{\max}
\end{equation}
\end{lemma}
Plugging Eq. \ref{eq:s_lower_bound} on the RHS of Eq. \ref{eq:cubic_model_bound} yields the sufficient decrease with the Newton step $\vect{s}_t$ which is used in the proof of the following main theorem of our algorithm.

% TODO seonho: first order optimality. From Lemma 4 change \epsilon to \epsilon^{1/2}. to attain \mathcal{O}(\epsilon^{-1.5} which is the state of the art...
%  TODO seonho: put some remark to mention that there is an opportunity to come up with \mathcal{O}\epsilon^{-1.75} for time complexity which is the state of the art.

\begin{theorem}[Worst-Case Iteration Complexity for Approximate First- and Second-Order Optimality]\label{theorem:worst_case_iteration_complexity}
Let A\ref{assump:function_boundness}, A\ref{assump:Lipschitz}, A\ref{assump:grad_hess_bound}, A\ref{assump:lanczos_termination}, and A\ref{assump:convergence} hold. 
For 
\begin{equation}\label{eq:epsilon_condition_custom}
    1 > \epsilon > \max \{3\epsilon_\vect{g},144\epsilon_\vect{B}^2,\frac{2}{\kappa_l}\epsilon_\vect{B}^2\},
\end{equation}
Algorithm \ref{alg:SANC} provides an iteration $t^*$ such that $\norm{\nabla f(\vect{x}_{t^*})} \leq \epsilon$ within at most $\mathcal{O}(\epsilon^{-3/2})$ iterations.

Also it provides an iteration $t^*$ such that $\norm{\nabla f(\vect{x}_{t^*})} \leq \epsilon$ and $-\lambda_{\min}(\vect{Q}_{t^*}^T\vect{B}_{t^*}\vect{Q}_{t^*}) \leq \epsilon$ within at most
\begin{equation}
    l_1+l_2+1=\mathcal{O}(\epsilon^{-3})
\end{equation}
iterations.
Both are with high probability $1-\delta$.
Here,
\begin{align}
    l_1:=\ceil[\bigg]{\frac{f(\vect{x}_0)-f_{low}}{\min\{\kappa_{s1}\epsilon^{\frac{3}{2}},\kappa_d\epsilon^3\}}},\label{eq:l1_def}\\
    l_2:=\ceil[\bigg]{\frac{f(\vect{x}_0)-f_{low}}{\min\{\kappa_{s2}\epsilon^3,\kappa_d\epsilon^3\}}},\label{eq:l2_def}\\
    \kappa_{s1}:=\frac{\eta_1\sigma_{\min}}{36\sqrt{6}\left( L_2+2\sigma_{\max} \right)^{\frac{3}{2}}},\,\,\:\:    \kappa_{s2}:=\frac{\eta_1\sigma_{\min}}{6\sigma_{\max}},\label{eq:pakkas_def}\\
    \text{and }\kappa_d:=\frac{1}{24L_2^2}\label{eq:kappad_def} %\kappa_u:=\frac{1}{\log\gamma_1}\log\left( \frac{\sigma_{\max}}{\sigma_{\min}} \right)\label{eq:kappau_def}
\end{align}
\end{theorem}

% \begin{theorem}[Worst-Case Iteration Complexity for Approximate First- and Second-Order Optimality]\label{theorem:worst_case_iteration_complexity}
% Let A\ref{assump:function_boundness}, A\ref{assump:Lipschitz}, A\ref{assump:grad_hess_bound}, A\ref{assump:lanczos_termination}, and A\ref{assump:convergence} hold. 
% Algorithm \ref{alg:SANC} provides an iteration $t^*$ such that $\norm{\nabla f(\vect{x}_{t^*})} \leq \epsilon$ and $-\lambda_{\min}(\vect{Q}_{t^*}^T\vect{B}_{t^*}\vect{Q}_{t^*}) \leq \epsilon$ within at most 
% \begin{equation}
%     l_1+l_2+1=\mathcal{O}(\epsilon^{-3})
% \end{equation}
% iterations with high probability $1-\delta$ for 
% \begin{equation}\label{eq:epsilon_condition_custom}
%     1 > \epsilon > \max \{3\epsilon_\vect{g},12\epsilon_\vect{B},\frac{2}{\kappa_l}\epsilon_\vect{B}^2\}
% \end{equation}
% where
% \begin{align}
%     l_1:=\ceil[\bigg]{\frac{f(\vect{x}_0)-f_{low}}{\min\{\kappa_{s1}\epsilon^{\frac{3}{2}},\kappa_d\epsilon^3\}}},\label{eq:l1_def}\\
%     l_2:=\ceil[\bigg]{\frac{f(\vect{x}_0)-f_{low}}{\min\{\kappa_{s2}\epsilon^3,\kappa_d\epsilon^3\}}},\label{eq:l2_def}\\
%     \kappa_{s1}:=\frac{\eta_1\sigma_{\min}}{36\sqrt{6}\left( L_2+2\sigma_{\max} \right)^{\frac{3}{2}}},\,\,\:\:    \kappa_{s2}:=\frac{\eta_1\sigma_{\min}}{6\sigma_{\max}},\label{eq:pakkas_def}\\
%     \text{and }\kappa_d:=\frac{1}{24L_2^2}\label{eq:kappad_def} %\kappa_u:=\frac{1}{\log\gamma_1}\log\left( \frac{\sigma_{\max}}{\sigma_{\min}} \right)\label{eq:kappau_def}
% \end{align}
% \end{theorem}

\textbf{Remark 1} It is noted that Eq \ref{eq:epsilon_condition_custom} is intuitively reasonable because gradient and Hessian approximation bounds assumption \ref{assump:grad_hess_bound} does not depend on the norm of $\vect{s}_t$.
Using smaller $\epsilon_\vect{g}$ and $\epsilon_\vect{B}$, iterates of the SANC algorithm can converge to a closer point to an exact optimum.

\textbf{Remark 2} Theorem \ref{theorem:worst_case_iteration_complexity} imply that the time complexity of the SANC algorithm is about $\mathcal{O}(\epsilon^{-3.5})$ when equipped with the Lanczos method for attaining steps, $\vect{s}_t$ or $\vect{d}_t$ to find a point satisfying second order optimality. 

\textbf{Remark 3} For the first order optimality, the iteration complexity $\mathcal{O}(\epsilon^{-1.5})$ of Algorithm \ref{alg:SANC} is same as that of ARC algorithm.
If it is equipped with a method with per iteration complexity of $\mathcal{O}(\epsilon^{-0.25})$, the SANC find first order optimality in a time $\mathcal{O}(\epsilon^{-1.75})$ as \cite{agarwal2017finding,carmon2018accelerated}.

All the proofs for the above theorems are provided in the appendix.
% However, we need to research further on the time complexity to make it precise.

\section{Empirical Studies}\label{sec:empirical_studies}

\subsection{Practical implementation of SANC}
As a default setting, we used $\gamma=2$, $\eta_1=0.2$, $\eta_2=0.8$, and $\sigma_0=1.$. 
These parameter settings are same as those in SCR \cite{kohler2017sub}.
For the neural networks problems, we used $\eta_1=0.1$, and $\eta_2=0.3$ tuned by some experiments.

$L_1$ and $L_2$ parameters in SANC have tuned.
The searching range for $L_1$ and $L_2$ are $10^{-3:1:3}$.
The used parameter values are presented in Table \ref{table:datasets}.

\begin{table*}[!ht]
\caption{Datasets used for numerical experiments. $n$ stands for the number of datapoints, $f$ stands for the number of features, and $c$ stands for the number of classes. $L_1$ and $L_2$ represent the parameters used in the SANC and NCD method.}
\label{table:datasets}
\vskip 0.15in
\begin{center}
\begin{tabular*}{0.8\textwidth}{@{}l|llllll@{}}
\toprule
        & $n$      & $f$                               & $c$ & learning model      & $L_1$  & $L_2$        \\ \midrule
w1a     & 2477     & 300                                             & 2   & logistic regression  & 10.0 &  10.0      \\
w8a     & 49749    & 300                                             & 2   & logistic regression  & 10.0 &  10.0      \\
a9a     & 32561    & 123                                             & 2   & logistic regression  & 10.0 &  10.0      \\
ijcnn1  & 49990    & 22                                              & 2   & logistic regression  & 10.0 &  100.0      \\
covtype & 581012   & 54                                              & 2   & logistic regression  & 10.0 &  100.0      \\
higgs   & 11000000 & 28                                              & 2   & logistic regression  & 10.0 &  100.0      \\
segment & 2310     & 19                                              & 3   & multi-layer neural networks & 100.0 & 100.0\\
seismic & 78823    & 50                                              & 7   & multi-layer neural networks & 100.0 & 100.0\\
MNIST   & 60000    & 28$\times$28$\times$1 & 10  & convolutional neural networks & 100.0 & 100.0  \\
CIFAR10 & 50000    & 32$\times$32$\times$3 & 10  & convolutional neural networks & 1000.0 & 1000.0 \\ \bottomrule
\end{tabular*}
\end{center}
\vskip -0.1in
\end{table*}

For the Lanczos method, we truncated it at Lanczos iteration $5$ for all experiments.
Also, an approximate minimizer of an approximate local cubic model over a Krylov subspace was obtained by the 'conjugate gradient method' in \texttt{scipy.optimize.minimize}.
To calculate the left-most eigenpair of $\vect{T}$, we used \texttt{eigh\textunderscore{}tridiagonal} in \texttt{scipy.linalg}.

For SANC and SCR, we have to calculate $f(\vect{x}_t)$ and $f(\vect{x}_t+\vect{s}_t)$ to measure $\rho_t$.
For the neural networks problems, because of the shortage of the memory storage, we made another set of dataset whose size is same as $\mathcal{S}_\vect{g}$ and $\mathcal{S}_\vect{B}$ to estimate $f(\vect{x}_t)$ and $f(\vect{x}_t+\vect{s}_t)$.
% TODO: set of dataset whose size is --> datasets whose size are ?

\subsection{Setup}\label{subsec:setup}
In this section, we present some experimental results to show the effectiveness of the SANC algorithm for stochastic nonconvex optimization problems.
In our numerical experiments, we considered three machine learning problems, namely,
\begin{enumerate*}[label=(\roman*)]
    \item the logistic regression
    \item the multi-layer neural networks and
    \item the convolutional neural networks (CNN)
\end{enumerate*}
with real datasets.

For the logistic regression problems, with the binary classification dataset, i.e., $\{\vect{x}_i,y_i\}$, $y_i\in\{0,1\}$, in order to find the optimum $\vect{w}^*$, we solved the following problem,
\begin{align}
    \min_{\vect{w}\in\mathbb{R}^d}&-\frac{1}{n}\sum_{i=1}^{n}y_i\log \left(\frac{1}{1+e^{-\vect{w}^T\vect{x}_i}}\right)\\
    &+\log\left(\frac{e^{\vect{w}^T\vect{x}_i}}{1+e^{-\vect{w}^T\vect{x}_i}}\right)+\lambda\Omega(\vect{w})
\end{align}
where $\Omega(\vect{w})=\sum_{i=1}^d \frac{w_i^2}{1+w_i^2}$, and $\lambda$ is a fixed regularization coefficient. 
We initialized all variables by $\vect{w}_0=1$.

For multi-layer neural networks, we used two hidden layers. 
The number of neurons in the first hidden layers is $300$ and the number of neurons in the second hidden layers is $500$.
At the outer layer, softmax functions were used. 
Cross entropy loss was used for obtaining an objective function.
We also added $l2$ norm as a convex regularization term with a coefficient $\lambda=0.01$.
Also, we used hyperbolic tangent functions as nonlinear activation functions.

For the CNN, we used two convolutional receptive filters.
The first one is $5\times5\times32$ dimensional and the second one is $5\times5\times64$ dimensional.
The fully connected layer was added at the output of the convolutional layers, which has $1000$ neurons.
Also, we used hyperbolic tangent functions for nonlinear activation functions.
We also added $l2$ norm as a convex regularization term with a coefficient $\lambda=0.01$.
% TODO: 위의 문단이랑 반복되므로 
% The settings for nonlinear activation functions and $l2$ norm are same with those of multi-layer neural networks.

All the variables including the weights and the biased in the multi-layer neural networks and the CNN were initialized by the Xavier initialization \cite{glorot2010understanding}.

% For the first experiment, the logistic regression with a nonconvex regularization was considered.
% For the second and third experiments, we considered the neural network problems. 
% For multi-layer neural networks, we used two hidden layers with hyperbolic tangent activation functions for solving multi-classification problems.
% For CNN, we used two consecutive convolutional filters along with max-pooling layers for image classification problems.

We used real datasets from libsvm \cite{chang2011libsvm} for the logistic regression problem and multi-layer neural networks.
For CNN, we used MNIST \cite{lecun1998gradient} and CIFAR10 datasets \cite{krizhevsky2009learning}.

The sizes of $\mathcal{S}_\vect{g}$ and $\mathcal{S}_\vect{B}$ are $\ceil{\text{the number of datapoints}/20}$ for the logistic regression problems.
For the multi-layer neural networks and CNN, we computed stochastic gradient and Hessian estimates using independently drawn mini-batches of size 128.
These batch size settings are used for SANC as well as all the baselines which are described below.

We compared our method, SANC, with various optimization methods described in Section \ref{sec:related_works}.
These include stochastic gradient descent(SGD), Sub-sampled Cubic Regularization (SCR) \cite{kohler2017sub}, Cubic Regularization (CR) \cite{nesterov2006cubic}, and Negative Curvature Direction (NCD) method \cite{liu2018adaptive}.
It is noted that the only difference between SANC and SCR is to use a direction of negative curvature. 
The details on the setting of the each algorithm is as follows,
\begin{itemize}
\item 
SGD: Fixed step-size is used.
We used 0.01 for the logistic regression problems and 0.001 for the multi-layer neural networks and CNN problems, which is the best choice among $10^{-3:1:3}$, respectively.

\item 
SCR \cite{kohler2017sub}: We used the same parameter values what we used for SANC. We used same sets, $\mathcal{S}_\vect{g}$, $\mathcal{S}_\vect{B}$ in SANC, which is different with the experiments in \cite{kohler2017sub} where they used increasing sized sets.
% TODO: what we used for SANC --> of SANC

\item
CR \cite{nesterov2006cubic}: For the fixed cubic coefficient $\sigma$, we used $5.$ for all problems.

\item
NCD \cite{liu2018adaptive}: We used the same parameter values what we used for SANC.
% TODO: what we used for SANC --> of SANC
% We did not consider any other applications on NCD.
\end{itemize}
All baselines are our implementations using \texttt{Tensorflow}.
% TODO: All baselines of our implementations are based on tensorflow.
% Please refer to the appendix for more detailed descriptions of the method and problem settings.

\subsection{Results}\label{subsec:results}

\begin{figure}[h!]
\centering
\captionsetup{justification=centering}
  \subcaptionbox{w1a dataset\\($n=2477$, $d=300$)}{\includegraphics[width = 1.55in]{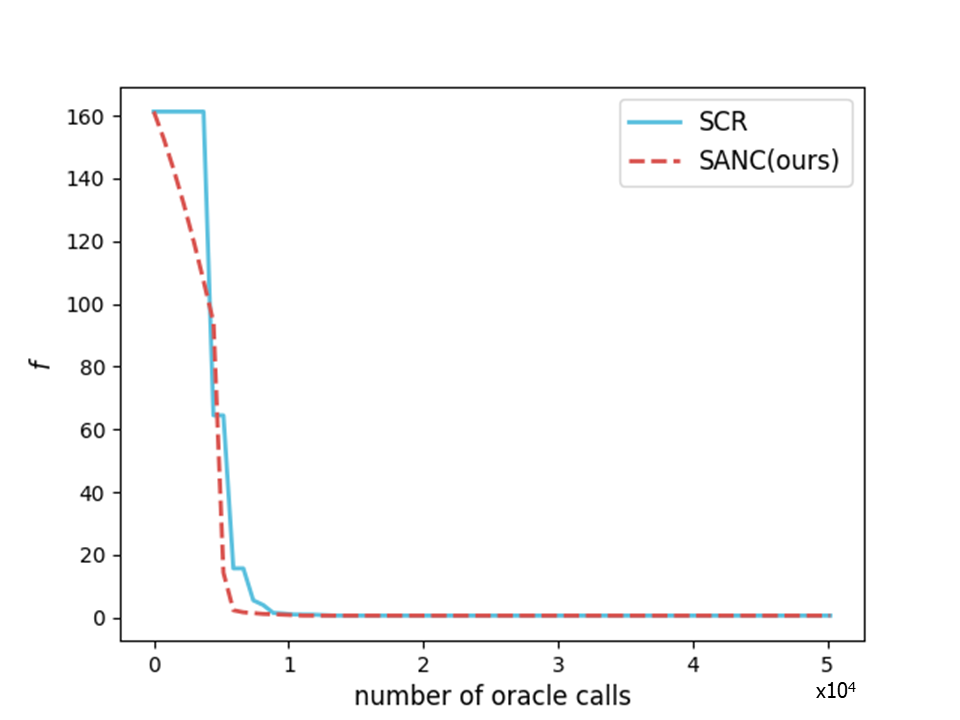}}\quad
  \subcaptionbox{higgs dataset\\($n=11000000$, $d=28$)}{\includegraphics[width = 1.55in]{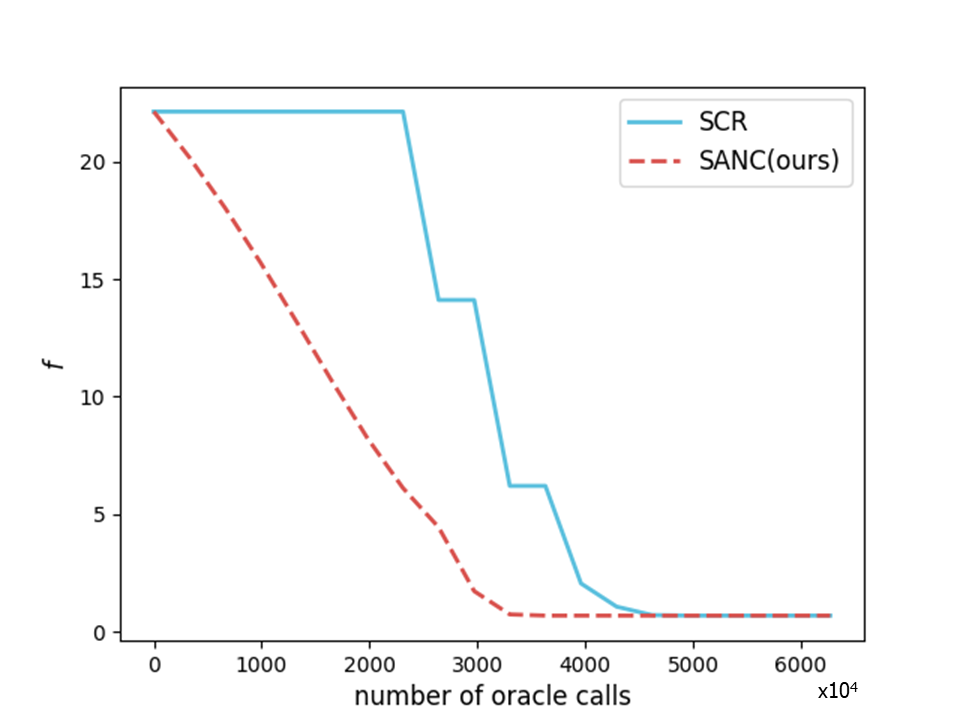}}\\
  \subcaptionbox{ijcnn1 dataset\\($n=49990$, $d=22$)}{\includegraphics[width = 1.55in]{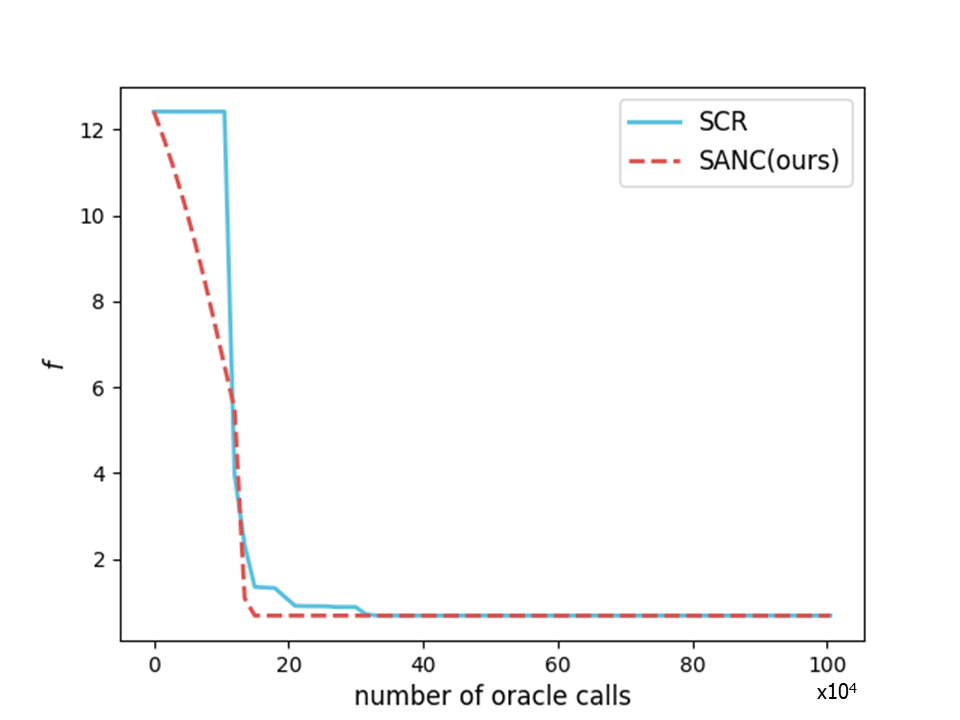}}\quad
  \subcaptionbox{covtype dataset\\($n=581012$, $d=54$)}{\includegraphics[width = 1.55in]{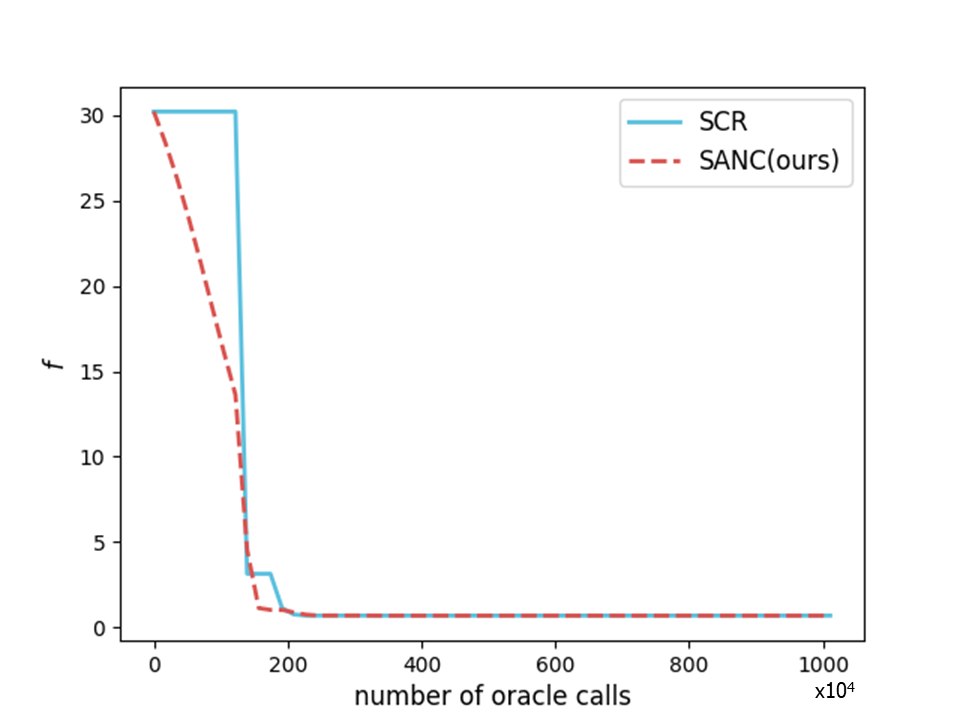}}
  \caption{Training losses of the logistic regression with a nonconvex regularization term with $\lambda=1.0$ over the number of oracle calls. We set $\sigma_0=0.001$ for both methods. This setting makes unsuccessful iterations at initial iterations. All function losses are the average of independent 10 runs.}\label{fig:lr_small_parameter_loss values}
\end{figure}

\begin{figure*}[!b]
\begin{subfigure}{0.333\textwidth}
% \centering
\includegraphics[height=3.7in]{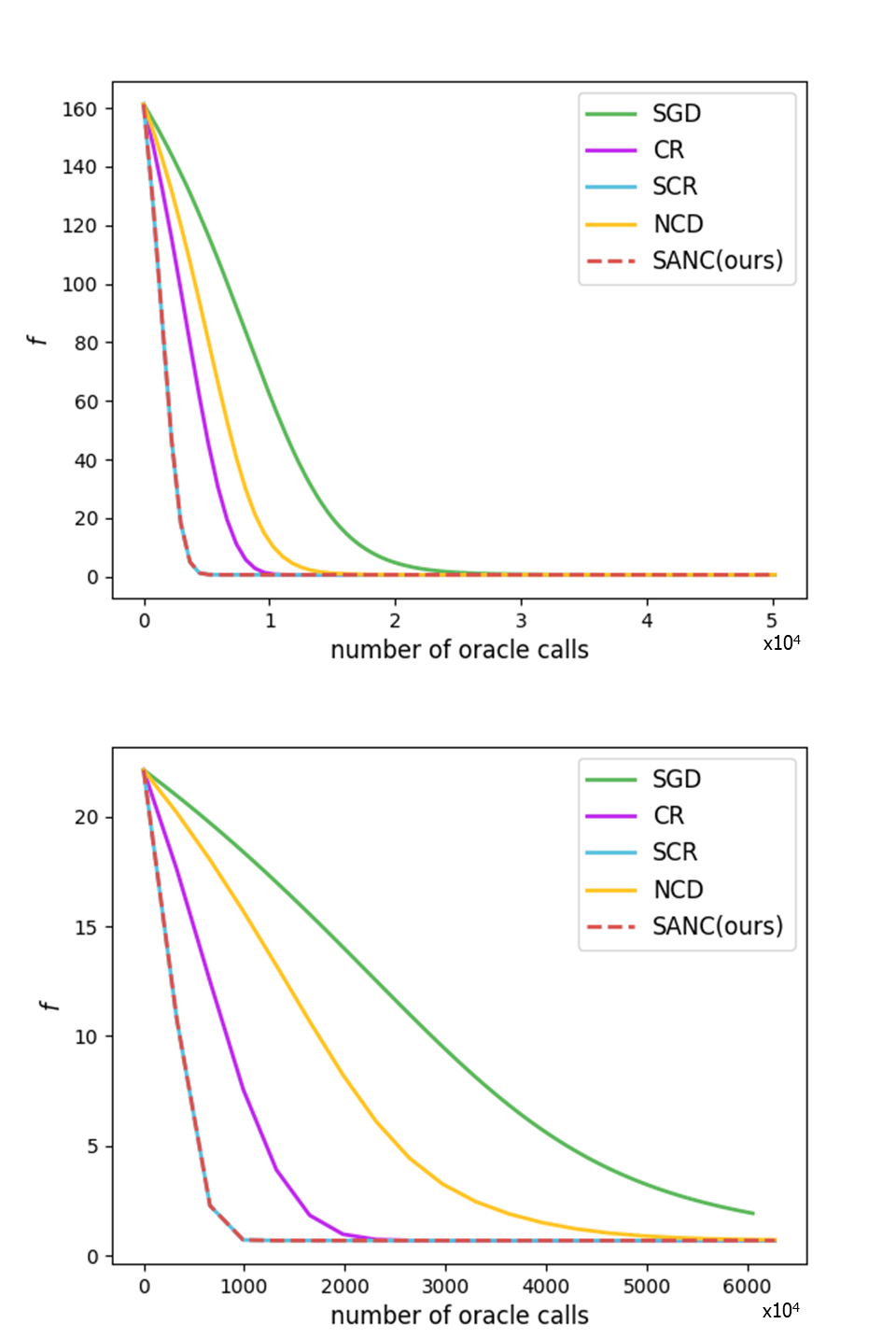}
\caption{Logistic regression}\label{subfig:lr}
\end{subfigure}\hfill
\begin{subfigure}{0.333\textwidth}
% \centering
\includegraphics[height=3.7in]{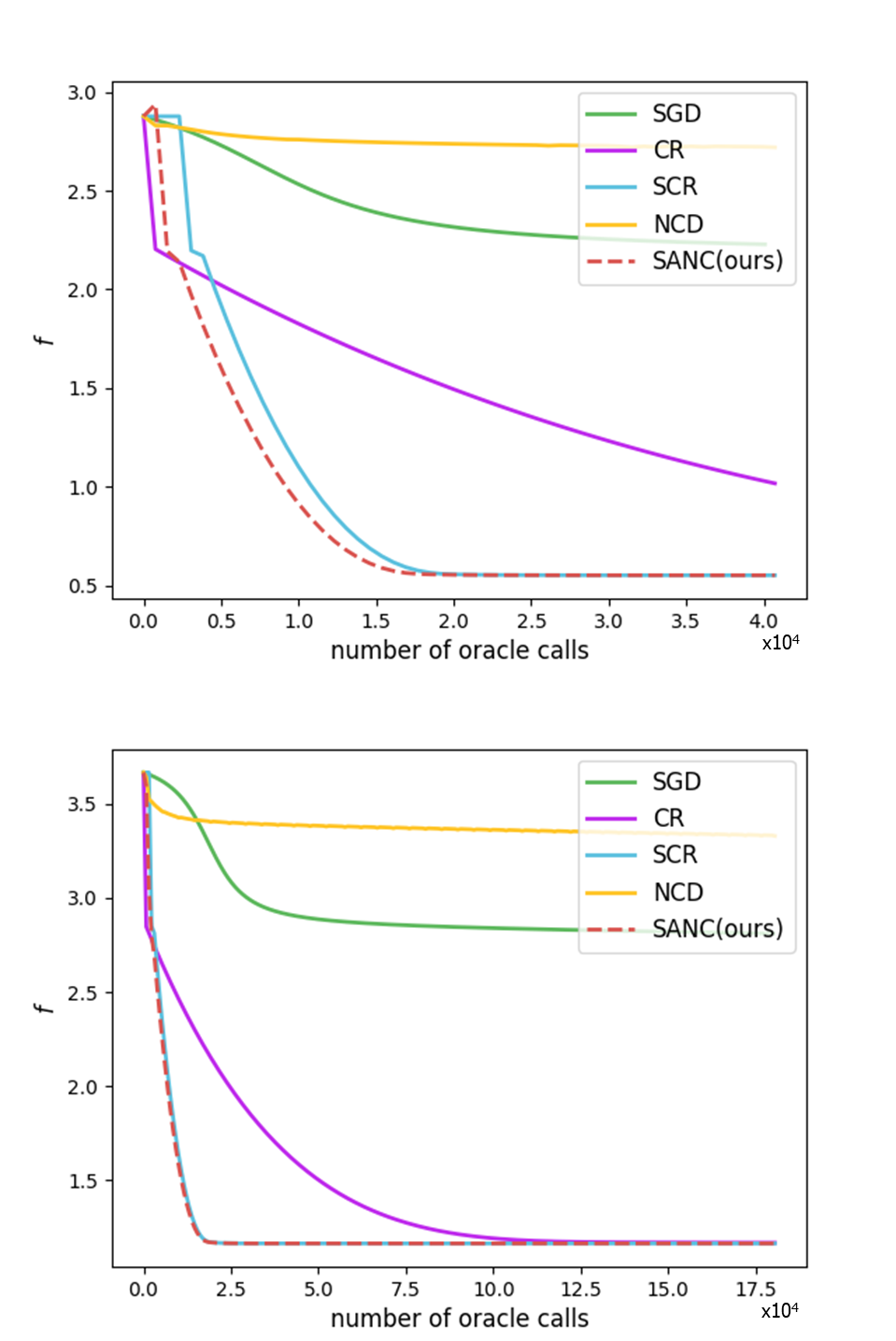}
\caption{Multi-layer neural networks}\label{subfig:mnn}
\end{subfigure}\hfill
\begin{subfigure}{0.333\textwidth}
% \centering
\includegraphics[height=3.7in]{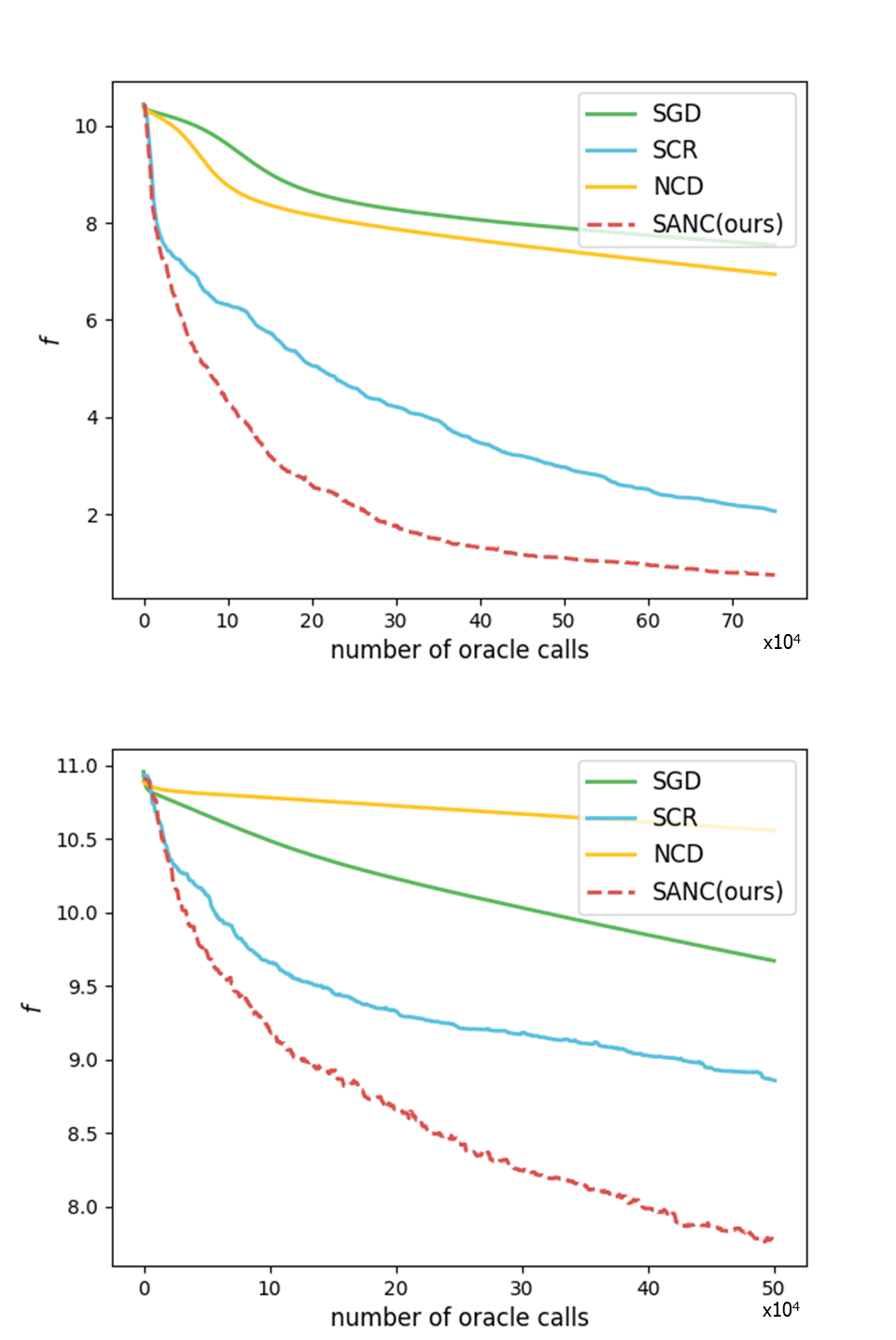}
\caption{Convolutional neural networks (CNN)}\label{subfig:cnn}
\end{subfigure}
\caption{Training losses over the number of oracle calls. (a) Logistic regression with a nonconvex regularization term with $\lambda=1.0$ solving for w1a ($n=2477$, $d=300$)(\textit{top}) and higgs ($n=11000000$, $d=28$)(\textit{bottom}) data. (b) multi-layer perceptron solving for seismic(\textit{top}) and segment(\textit{bottom}) data. (c) convolutional neural networks solving for MNIST(\textit{top}) and CIFAR10(\textit{bottom}) data. We omit CR method for CNN problems because of its poor performance. Please refer to the supplementary material for the details on the parameter and problem settings. All function losses are the average of independent 10 runs.
}\label{fig:loss values}
\end{figure*}

Figure \ref{fig:lr_small_parameter_loss values} shows that the SANC method has advantages to use a direction of negative curvature at unsuccessful iterations.
Because of the existence of the nonconvex regularization term, initial iterates have indefinite Hessians.
As pointed out in Section \ref{sec:SANC}, with a smaller $\sigma_0$, we experienced that the SCR method suffers from unsuccessful iterations in initial iterations. 

Figure \ref{fig:loss values} illustrates the training loss results for three problems.
It shows that SANC and SCR are identical for the logistic regression problems. 
This is because SANC does not make the update with a direction of negative curvature in that problems.
As shown in Figure \ref{subfig:mnn}, \ref{subfig:cnn}, SANC reveals competitive loss history against others.
Because neural networks problems have so many local minima and saddle points, it is suitable to use a negative curvature step to find a better step at unsuccessful iteration.
Please refer to the supplementary material for more numerical results.

\section{Conclusion}\label{sec:conclusion}
In this work, we have developed a stochastic optimization framework which is called SANC for nonconvex optimization. 
The novelty of the proposed algorithms is that we combine the adaptive cubic regularized Newton method and the negative curvature method.
This method not only maintain the trust-region like framework of the adaptive cubic regularized Newton method but also could reduce the number of oracle calls for the large-scale machine learning problems with real datasets.
We would like to point out that a Krylov subspace where an approximate minimizer of an approximate local cubic model is established can also be efficiently used for approximating the left-most eigenvector of the Hessian.
This eigenvector can be utilized as a direction of negative curvature.
Also, as shown in empirical results, a subsampling technique, that we proposed, is suitable to be used for solving large-scale machine learning problems.
As far as we know, this approach is the first work to integrate the negative curvature method into the Newton method. 
We think the accelerated method and the stochastic variance reduced gradient method can perform well with our algorithm to attain better theoretical and empirical results.

\section*{acknowledgements}
% Acknowledgements to sponsoring agencies and individuals should be placed here.
%TODO: 
This work was partially supported by the National Research Foundation (NRF) of Korea. (NRF-2018R1D1A1B07043406).
% \end{acknowledgements}

% TODO: (Acknowledgement) This work was supported by the National Research Foundation (NRF) of Korea. (NRF-2018R1D1A1B07043406). 

\bibliography{ref}
\bibliographystyle{icml2019}

%%%%%%%%%%%%%%%%%%%%%%%%%%%%%%%%%%%%%%%%%%%%%%%%%%%%%%%%%%%%%%%%%%%%%%%%%%%%%%%
%%%%%%%%%%%%%%%%%%%%%%%%%%%%%%%%%%%%%%%%%%%%%%%%%%%%%%%%%%%%%%%%%%%%%%%%%%%%%%%
\newpage
\appendix
\onecolumn
\section{Appendix}
For the completeness, we rewrite all the lemmas and the theorems of the article.

\subsection{Gradient and Hessian sampling bound}
Before addressing Theorem \ref{theorem:gradient_sampling_bound} in the article, we introduce the vector Bernstein inequality and its following lemma. These proofs are based on the work from \cite{kohler2017sub}. 
\begin{customlemma}{8}[Vector Bernstein Inequality, Theorem 12 in \cite{gross2011recovering}]\label{lem:vector_bernstein_inequality}
Let $\vect{X}_1,\dots\vect{X}_m$ be independent zero-mean vector-valued random variables. Let
\begin{equation}
    N=\norm{\sum_{i=1}^m\vect{X_i}}_2
\end{equation}
then
\begin{equation}\label{eq:vector_bernstein}
    \Pr[N\geq\sqrt{V}+t] \leq \exp{\left( -\frac{t^2}{4V} \right)}
\end{equation}
where $V=\sum_i \mathbb{E}\left[ \norm{\vect{X}_i}_2^2 \right]$ and $t \leq V/max_i\norm{\vect{X}_i}_2$.
\end{customlemma}
\begin{proof}
please refer to the proof of the theorem 12 in \cite{gross2011recovering}
\end{proof}

Before describing the proof for Theorem \ref{theorem:gradient_sampling_bound} in the paper, we describe the gradient bound in terms of the sampling size to establish $\epsilon_\vect{g}$ error.

\begin{customlemma}{9}[Approximate Gradient Deviation Bound]\label{lem:grad_deviation}
For $\epsilon_\vect{g} \leq 2L_0\left(1+\frac{1}{\sqrt{|\mathcal{S}_{\vect{g}}|}} \right)$, if 
\begin{equation}
    \norm{\vect{g}_t-\nabla f(\vect{x}_t)} \leq 2L_0\frac{1+2\sqrt{\log{\frac{1}{\delta}}}}{\sqrt{|\mathcal{S}_{\vect{g}}|}}
\end{equation}
then for small $\epsilon_\vect{g}>0$ we can establish the condition \ref{eq:grad_bound}
with high probability $1-\delta$
\end{customlemma}
\begin{proof}
We use Lemma \ref{lem:vector_bernstein_inequality} for the vector Bernstein inequality.
Before using it, first, we make independent zero-mean vector-valued random variables,
\begin{equation}
    \vect{X}_i=\nabla f_i(\vect{x})-\nabla f(\vect{x})
\end{equation}
Then, 
\begin{equation} \label{eq:norm_zero_mean_vector}
    \norm{\vect{X}_i}=\norm{\nabla f_i(\vect{x})-\nabla f(\vect{x})}\leq \norm{\nabla f_i(\vect{x})}+\norm{\nabla f(\vect{x})} \leq 2L_0
\end{equation}
where the first inequality is from the triangular inequality and the second inequality is from the Lipschitz continuity of the function values.
Hence, 
\begin{equation}
    \mathbb{E}\left[ \norm{\vect{X}_i}^2 \right] \leq 4L_0^2
    \Rightarrow V=\sum_{i\in \mathcal{S}_\vect{g}}\mathbb{E}\left[ \norm{\vect{X}_i}^2 \right] \leq 4L_0^2|\mathcal{S}_{\vect{g}}|
\end{equation}
From Eq. \ref{eq:vector_bernstein}, let $|\mathcal{S}_{\vect{g}}|\epsilon_\vect{g}:=\sqrt{V}+t$ then,
\begin{align}
    \Pr\left[\norm{\sum_{i\in\mathcal{S}_\vect{g}}\left(\nabla f_i(\vect{x})-\nabla f(\vect{x}) \right)} \geq |\mathcal{S}_\vect{g}|\epsilon_\vect{g}\right] &\leq \exp\left({-\frac{1}{4}\left( \frac{|\mathcal{S}_\vect{g}|\epsilon_\vect{g}}{\sqrt{V}}-1 \right)^2}\right)\\
    &\leq \exp\left({-\frac{1}{4}\left( \frac{\sqrt{|\mathcal{S}_\vect{g}|}\epsilon_\vect{g}}{2L_0}-1 \right)^2 }\right) \leq \delta\\
    &\Leftrightarrow \frac{1}{4}\left( \frac{\sqrt{|\mathcal{S}_\vect{g}|}\epsilon_\vect{g}}{2L_0}-1 \right)^2\geq \log\left(\frac{1}{\delta}\right)\\
    &\Leftrightarrow \epsilon_\vect{g} \geq \frac{\left(1+\sqrt{4\log(1/\delta)}\right)2L_0}{\sqrt{|\mathcal{S}_\vect{g}}|}
\end{align}
Thus, the probability of $\epsilon_\vect{g} < \frac{\left(1+\sqrt{4\log(1/\delta)}\right)2L_0}{\sqrt{|\mathcal{S}_\vect{g}}|}$ is higher or equal to $1-\delta$.
Rearranging this condition yields the conclusion.
\end{proof}

\begin{customthm}{1}[Gradient Sampling Bound]
If
\begin{equation}
|\mathcal{S}_{\vect{g},t}| \geq  \frac{4L_0^2\left( 1+2\sqrt{\log{\frac{1}{\delta}}} \right)^2}{\epsilon_\vect{g}^2},
\end{equation}
then $\vect{g}_t$ satisfies the gradient approximation bound Eq. \ref{eq:grad_bound} for all $t\geq0$ with high probability $1-\delta$. 
\end{customthm}
\begin{proof}
From Lemma \ref{lem:grad_deviation}, 
\begin{align}
    \frac{\left(1+\sqrt{4\log(1/\delta)}\right)2L_0}{\sqrt{|\mathcal{S}_\vect{g}}|} \leq \epsilon_\vect{g}\\
    \Rightarrow \sqrt{|\mathcal{S}_\vect{g}|} \geq \frac{2L_0(1+\sqrt{4\log(1/\delta)}}{\epsilon_\vect{g}}.\\
\end{align}
Squaring both sides finalizes the proof.
\end{proof}

For the Hessian sampling size bound, we use a similar approach.
The proof is based on the Operator Bernstein Inequality. 
The big difference between the gradient sampling and Hessian sampling is that Hessian sampling size depends on the number of variables, $d$.

\begin{customlemma}{10}[Operator-Bernstein Inequality, Theorem 6 in \cite{gross2011recovering}]
Let $\vect{X}_i,\,i=1,\dots,m$ be i.i.d., $n$ by $n$, zero-mean, Hermitian matrix-valued random variables.
Assume $V_0,c\in\mathbb{R}$ are such that $\norm{\mathbb{E}[\vect{X}_i^2]}\leq V_0^2$ and $\norm{\vect{X}_i}\leq c$. 
Set $S=\sum_{i=1}^m \vect{X}_i$ and let $V=mV_0^2$.
Then,
\begin{equation}\label{eq:operator_bernstein_inequality}
    \Pr[\norm{S}>t] \leq 2n \exp\left( -\frac{t^2}{4V} \right)
\end{equation}
for $t\leq 2V/c$.
\end{customlemma}
\begin{proof}
Please refer to the proof of the theorem 6 in \cite{gross2011recovering}.
\end{proof}

Let us denote $\vect{H}_i$ as an abbreviation of $\nabla^2 f_i(\vect{x})$ for some $i$ and $\vect{H}$ denote the expectation of $\vect{H}_i$.
\begin{customlemma}{11}[Approximate Hessian Deviation Bound]\label{lem:hessian_deviation_bound}
For $\epsilon_\vect{B} \leq 4L_1$, if
\begin{equation}
    \norm{\vect{B}-\vect{H}} < \frac{4L_1\sqrt{\log(2d/\delta)}}{\sqrt{|\mathcal{S}_\vect{B}|}}
\end{equation}
then for small $\epsilon_\vect{B}>0$ we can establish the condition \ref{eq:hessian_bound}
with high probability $1-\delta$.
\end{customlemma}
\begin{proof}
To use the operator-Bernstein inequality \ref{eq:operator_bernstein_inequality}, We introduce an independent zero mean Hermitian random variable,
\begin{equation}
    \vect{X}_i=\vect{H}_i-\vect{H}
\end{equation}
Then,
\begin{equation}
    \norm{\vect{X}_i} = \norm{\vect{H}_i-\vect{H}} \leq 2L_1
\end{equation}
where the inequality is from the triangular inequality and the fact that $f_i$ has Lipschitz continuous gradient.
Let define $V_0^2:=4L_1^2$ and $V:=|\mathcal{S}_\vect{B}|V_0^2=|\mathcal{S}_\vect{B}|4L_1^2$
Then,
\begin{equation}
    \norm{\vect{B}-\vect{H}} = \frac{1}{|\mathcal{S}_\vect{B}|}\norm{\sum_{i\in\mathcal{S}_\vect{B}}{\vect{X}_i}}
\end{equation}
Hence from the operator-Bernstein inequality \ref{eq:operator_bernstein_inequality},
\begin{equation}
    \Pr\left[\norm{\vect{B}-\vect{H}} \geq \epsilon_\vect{B}\right]
    =\Pr\left[\norm{\sum_{i\in\mathcal{S}_\vect{B}}{\vect{X}_i}}\geq |\mathcal{S}_\vect{B}|\epsilon_\vect{B}\right] \leq 2d\cdot \exp\left( -\frac{\epsilon_\vect{B}^2|\mathcal{S}_\vect{B}|}{16L_1^2} \right) \leq \delta
\end{equation}
for $|\mathcal{S}_\vect{B}|\epsilon_\vect{B}\leq \frac{4|\mathcal{S}_\vect{B}|L_1^2}{L_1}$ that means $\epsilon_\vect{B}\leq 4L_1$.
Thus,
\begin{align}
    2d\cdot \exp\left(-\frac{\epsilon_\vect{B}|\mathcal{S}_\vect{B}|}{16L_1^2}  \right) \leq \delta\\
    \Leftrightarrow \frac{\epsilon_\vect{B}^2|\mathcal{S}_\vect{B}|}{16L_1^2} \geq \frac{2d}{\delta}\\
    \Leftrightarrow \epsilon_\vect{B} \geq \frac{4L_1\sqrt{\log(2d/\delta)}}{\sqrt{|\mathcal{S}_\vect{B}|}}
\end{align}
\end{proof}

From this result, we can prove the Hessian sampling bound.

\begin{customthm}{2}[Hessian Sampling Bound]
If
\begin{equation}
    |\mathcal{S}_{\vect{B},t}| \geq \frac{16L_1^2\log{\frac{2d}{\delta}}}{\epsilon_\vect{B}^2},
\end{equation}
then $\vect{B}_t$ satisfies the Hessian approximation bound Eq. \ref{eq:hessian_bound} for all $t\geq0$ with high probability $1-\delta$.
\end{customthm}
\begin{proof}
From the previous lemma \ref{lem:hessian_deviation_bound}, 
\begin{align}
    4L_1\sqrt{\frac{\log(2d/\delta)}{|\mathcal{S}_\vect{B}|}} \leq \epsilon_\vect{B}\\
    \Leftrightarrow \frac{\log(2d/\delta)}{|\mathcal{S}_\vect{B}|} \leq \left( \frac{\epsilon_\vect{B}}{4L_1} \right)^2\\
    \Leftrightarrow |\mathcal{S}_\vect{B}| \geq \frac{16L_1^2\log{\frac{2d}{\delta}}}{\epsilon_\vect{B}^2}
\end{align}
\end{proof}

\subsection{Worst-case iteration complexity bound for achieving approximate first- and second-order optimality}\label{subsec:worst_case_complexity}

\begin{customlemma}{3}[Cubic Regularization Coefficient Bound]
Let the A\ref{assump:Lipschitz}, A\ref{assump:grad_hess_bound} and A\ref{assump:convergence} hold. 
For finite positive values $\sigma_{\min}$ and $\sigma_{\max}$,
it holds,
\begin{equation}
    \sigma_{\min} \leq \sigma_t \leq \sigma_{\max} \: \text{ for all }t\geq0
\end{equation}
where $\sigma_{\max}:=\max\left\{\sigma_0,\gamma\left( \frac{3}{2}L_2+\frac{3(\epsilon_\vect{g}+\frac{1}{2}\epsilon_\vect{B})}{\epsilon_s^2} \right) \right\}$.
\end{customlemma}
\begin{proof}
For $\sigma_{\min}$, $\sigma_t$ is always bounded below by $\epsilon_m$, the machine precision.
Further, we need to show that $\sigma_t\leq \sigma_{\max}$.
For each iteration $t\geq0$,
\begin{align}
f(\vect{x}_t+\vect{s}_t)-\widetilde{m}_t(\vect{s}_t)&=(\nabla f(\vect{x}_t)-\vect{g}_t)^T\vect{s}_t+\frac{1}{2}\vect{s}_t^T(\vect{H}(\xi_t)-\vect{B}_t)\vect{s}_t-\frac{1}{3}\sigma_t\norm{\vect{s}_t}^3\\
&\leq \epsilon_\vect{g}\norm{\vect{s}_t}+\frac{1}{2}\norm{\vect{s}_t}^2\norm{\vect{H}(\xi_t)-\vect{H}(\vect{x}_t)}+\frac{1}{2}\norm{\vect{H}(\vect{x}_t)-\vect{B}_t}\norm{\vect{s}_t}-\frac{1}{3}\sigma_t\norm{\vect{s}_t}^3\\
&\leq \epsilon_\vect{g}\norm{\vect{s}_t}+\frac{1}{2}L_2\norm{\vect{s}_t}^3+\frac{1}{2}\epsilon_\vect{B}\norm{\vect{s}_t}-\frac{1}{3}\sigma_t\norm{\vect{s}_t}^3\\
&=\left( \epsilon_\vect{g}+\frac{1}{2}\epsilon_\vect{B} \right)\norm{\vect{s}_t}+\left( \frac{1}{2}L_2-\frac{\sigma_t}{3} \right)\norm{\vect{s}_t}^3\label{eq:last_coeff_bound}
\end{align}
where $\xi_t$ is a line segment of $\vect{x}_t$ and $\vect{x}_t+\vect{s}_t$.

To achieve $f(\vect{x}_t+\vect{s}_t)\leq \widetilde{m}_t(\vect{s}_t)$ from Eq. \ref{eq:last_coeff_bound}, we need to have the condition that 
\begin{align}
\left( \epsilon_\vect{g}+\frac{1}{2}\epsilon_\vect{B} \right)\norm{\vect{s}_t}+\left( \frac{1}{2}L_2-\frac{\sigma_t}{3} \right)\norm{\vect{s}_t}^3 \leq 0\\
% \Rightarrow \sigma_t \geq \frac{3}{2}L_2+\frac{3}{\norm{\vect{s}_t}^2}\left( \epsilon_\vect{g}+\frac{1}{2}\epsilon_\vect{B} \right) \geq \frac{3}{2}L_2+\frac{3}{\epsilon_s^2}\left( \epsilon_\vect{g}+\frac{1}{2}\epsilon_\vect{B} \right)
\Rightarrow \sigma_t \geq  \frac{3}{2}L_2+\frac{3}{\epsilon_s^2}\left( \epsilon_\vect{g}+\frac{1}{2}\epsilon_\vect{B} \right)
\end{align}
% $\sigma_t \geq \frac{3}{2}L_2+\frac{3(\epsilon_\vect{g}+\frac{1}{2}\epsilon_\vect{B})}{\epsilon_s^2}$.
Then, from the definition of $\rho_t$, the iteration $t$ is very successful and $\sigma_{t+1}\leq \sigma_t$. 
Let us think that $\sigma_t$ is slightly less than $\frac{3}{2}L_2+\frac{3(\epsilon_\vect{g}+\frac{1}{2}\epsilon_\vect{B})}{\epsilon_s^2}$, then it could be multiplied by $\gamma$ and the case iterations can be started with a large $\sigma_0$. 
This gives the desired result.
\end{proof}

For the two consecutive iterations, we have to establish a sufficient decrease.
Hence we provide the following lemma for the negative curvature step.
\begin{customlemma}{4}[Sufficient Decrease with Negative Curvature Step, Lemma 3 in \cite{liu2018adaptive}]
When $\vect{v}_t^T\vect{B}_t\vect{v}_t\leq -\epsilon'/2$ and $\epsilon_\vect{B} \leq \epsilon'/12$, the negative curvature step $\vect{d}_t$ satisfies that 
\begin{equation}\label{eq:d_update_bound_apdx}
    f(\vect{x}_t)-\mathbb{E}[f(\vect{x}_t+\vect{d}_t)] \geq \max\left\{ \frac{{\epsilon'}^3}{24L_2^2},\frac{\norm{\vect{g}_t}^2}{4L_1}-\frac{\epsilon_\vect{g}^2}{L_1} \right\}
\end{equation}
with high probability $1-\delta$.
\end{customlemma}

\begin{proof}
Please refer to the proof of the Lemma 3 in \cite{liu2018adaptive}.
\end{proof}

From the paper, we showed that the approximate minimizer $s_t$ over the Krylov subspace $\mathcal{K}_j$ satisfies Eq. \ref{eq:s_condition1} and Eq. \ref{eq:s_condition2}.
The next lemma gives a lower bound on the approximate local cubic model decrease when Eq. \ref{eq:s_condition1} and Eq. \ref{eq:s_condition2} are satisfied for a successful iteration $t\in\mathcal{N}$.
\begin{customlemma}{5}[Local Cubic Model Decrease, Lemma 3.3 in \cite{cartis2011adaptive}]
Suppose that $\vect{s}_t$ satisfies Eq. \ref{eq:s_condition1} and Eq. \ref{eq:s_condition2}. Then for all (very) successful iterations $t$, 
\begin{equation}
    f(\vect{x}_t)-\widetilde{m}_t(\vect{s}_t) \geq \frac{1}{6}\sigma_t\norm{\vect{s}_t}^3
\end{equation}
\end{customlemma}
\begin{proof}
Please refer to the proof of the Lemma 3.3 in \cite{cartis2011adaptive}. Please note that in \cite{cartis2011adaptive} they established an approximate local cubic model with exact gradient $\nabla f$ and approximate Hessian $\vect{B}$, but in our approach we establish an approximate local cubic model with approximate gradient $\vect{g}$ and Hessian $\vect{B}$. 
However, the procedure of the proof is same except for using $\vect{g}$ instead of $\nabla f$.
\end{proof}

Also, we establish the relationship between the magnitude of $\vect{s}_t$ and $\vect{g}_t$ for a successful iteration $t$. 

\begin{customlemma}{6}[Lower Bound of $\vect{s}_t$]
Suppose that A\ref{assump:Lipschitz}, A\ref{assump:grad_hess_bound}, A\ref{assump:lanczos_termination}, and A\ref{assump:convergence} hold. 
For all (very) successful iterations $t$, the Newton step $\vect{s}_t$ satisfies that
\begin{equation}
    \norm{\vect{s}_t} \geq \frac{1}{\kappa_l}    \left(-\epsilon_\vect{B}+\sqrt{\epsilon_\vect{B}^2+2\kappa_l\left( \norm{\nabla f(\vect{x}_t+\vect{s}_t)-2\epsilon_\vect{g}} \right)} \right)
\end{equation}
where 
\begin{equation}
    \kappa_l:=L_2+2\sigma_{\max}
\end{equation}
\end{customlemma}
\begin{proof}
For a successful iteration $t$,
\begin{align}
\norm{\nabla f(\vect{x}_t+\vect{s}_t)} &\leq \norm{\nabla f(\vect{x}_t+\vect{s}_t)-\nabla \widetilde{m}_t(\vect{s}_t)}+\norm{\nabla \widetilde{m}_t(\vect{s}_t)}\\
 &= \norm{\nabla f(\vect{x}_t) + \vect{H}(\xi_t)\vect{s}_t- \vect{g}_t-\vect{B}_t\vect{s}_t-\sigma_t\norm{\vect{s}_t}\vect{s}_t}+\norm{\nabla \widetilde{m}_t(\vect{s}_t)}\\
 & = \norm{\nabla f(\vect{x}_t)-\vect{g}_t +\left( \vect{H}(\xi_t)-\vect{H}(\vect{x}_t) \right)\vect{s}_t+\left(\vect{H}(\vect{x}_t)-\vect{B}_t\right)\vect{s}_t - \sigma_t\norm{\vect{s}_t}\vect{s}_t} + \norm{\nabla \widetilde{m}_t(\vect{s}_t)}\\
 & \leq \epsilon_\vect{g} + \frac{1}{2}L_2\norm{\vect{s}_t}^2+\epsilon_\vect{B}\norm{\vect{s}_t} +  \sigma_{\max}\norm{\vect{s}_t}^2+\norm{\nabla \widetilde{m}_t(\vect{s}_t)}\\
 &= \left( \frac{1}{2}L_2+\sigma_{\max} \right)\norm{\vect{s}_t}^2+\epsilon_\vect{g}+\epsilon_\vect{B}\norm{\vect{s}_t}+\norm{\nabla \widetilde{m}_t(\vect{s}_t)}
\end{align}
where $\xi_t$ is a line segment between $\vect{x}_t$ and $\vect{x}_t+\vect{s}_t$. First inequality is from triangular inequality and second inequality is from Assumption \ref{assump:grad_hess_bound} and Eq. \ref{eq:sigma_limit}.

Hence with the TC \ref{assump:lanczos_termination},
\begin{equation}\label{eq:s_t_bound_nabla}
    \norm{\vect{s}_t} \geq \frac{1}{L_2+2\sigma_{\max}}\left( -\epsilon_\vect{B}+\sqrt{\epsilon_\vect{B}^2+(2L_2+4\sigma_{\max})\left( \norm{\nabla f(\vect{x}_t+\vect{s}_t)}-2\epsilon_\vect{g} \right)} \right)
\end{equation}
\end{proof}

\begin{customlemma}{12}\label{lemma:arithmetic_ineq}
With real numbers $a>0$, $b>0$, $t\in (0,1)$, $-a+\sqrt{a^2+b^2} \geq tb$ holds if and only if $b\geq \frac{2t}{\sqrt{1-t^2}}a$.
\end{customlemma}
\begin{proof}
\begin{align}
    -a+\sqrt{a^2+b^2} \geq tb\\
    \Leftrightarrow 2a^2+(1-t^2)b^2 \geq 2a\sqrt{a^2+b^2}\\
    \Leftrightarrow 4a^4+(1-t^2)^2b^4+4(1-t^2)a^2b^2 \geq 4a^4+4a^2b^2\\
    \Leftrightarrow (1-t^2)b^2-4t^2a^2 \geq 0\\
    \Leftrightarrow b \geq \frac{2t}{\sqrt{1-t^2}}a
\end{align}
\end{proof}

\begin{customthm}{7}[Worst-Case Iteration Complexity for Approximate First- and Second-Order Optimality]
Let A\ref{assump:function_boundness}, A\ref{assump:Lipschitz}, A\ref{assump:grad_hess_bound}, and A\ref{assump:lanczos_termination} hold. 
Algorithm \ref{alg:SANC} provides an iteration $t^*$ such that $\norm{\nabla f(\vect{x}_{t^*})} \leq \epsilon$ and $-\lambda_{\min}(\vect{Q}_{t^*}^T\vect{B}_{t^*}\vect{Q}_{t^*}) \leq \epsilon$ within at most 
\begin{equation}
    % (1+\kappa_u)(l_1+l_2)+2\kappa_u+1=\mathcal{O}(\epsilon^{-3})
    l_1+l_2+1=\mathcal{O}(\epsilon^{-3})
\end{equation}
iterations with high probability $1-\delta$ for 
\begin{equation}\label{eq:epsilon_condition_custom_apdx}
    1 > \epsilon > \max \{3\epsilon_\vect{g},12\epsilon_\vect{B},\frac{2}{\kappa_l}\epsilon_\vect{B}^2\}
\end{equation}
where
\begin{align}
    l_1:=\ceil[\bigg]{\frac{f(\vect{x}_0)-f_{low}}{\min\{\kappa_{s1}\epsilon^{\frac{3}{2}},\kappa_d\epsilon^3\}}},\\
    l_2:=\ceil[\bigg]{\frac{f(\vect{x}_0)-f_{low}}{\min\{\kappa_{s2}\epsilon^3,\kappa_d\epsilon^3\}}},\\
    \kappa_{s1}:=\frac{\eta_1\sigma_{\min}}{36\sqrt{6}\left( L_2+2\sigma_{\max} \right)^{\frac{3}{2}}},\,\,\:\:    \kappa_{s2}:=\frac{\eta_1\sigma_{\min}}{6\sigma_{\max}}\epsilon^3,\\
    \text{and }\kappa_d:=\frac{1}{24L_2^2} %\kappa_u:=\frac{1}{\log\gamma_1}\log\left( \frac{\sigma_{\max}}{\sigma_{\min}} \right)\label{eq:kappau_def}
\end{align}
\end{customthm}

\begin{proof}
\begin{enumerate}[(a)]
\item \label{enum:firstorder}
First, we need to show that an first order stationary point exists within at most the certain number of iterations.

For a sufficiently large iteration $j$, suppose that $\norm{\nabla f(\vect{x}_{t+1})} > \epsilon, \,\forall t\leq j$.
% and define $\mathcal{M}_j=\mathcal{N}_j \cup \mathcal{D}_j$.\\
From lemma \ref{lemma:arithmetic_ineq}, plugging $\epsilon_\vect{B}$, $\sqrt{\frac{2L_2+4\sigma_{\max}}{3}\epsilon}$ and $1/2$ into $a$, $b$ and $t$, respectively, yields
\begin{equation}\label{eq:epsb_eps}
-\epsilon_\vect{B}+\sqrt{\epsilon_\vect{B}+ \frac{2L_2+4\sigma_{\max}}{3} \epsilon} \geq \frac{1}{2}\sqrt{\frac{2L_2+4\sigma_{\max}}{3}\epsilon} = \sqrt{\frac{L_2+2\sigma_{\max}}{6}}\epsilon^{\frac{1}{2}}    
\end{equation}
if and only if $\epsilon \leq \frac{2}{\kappa_l}\epsilon_{\vect{B}}$ holds.

From the definition of $\epsilon$ Eq. \ref{eq:epsilon_condition_custom_apdx}, we assume that $\norm{\nabla f(\vect{x}_t+\vect{s}_t)}>\epsilon$, then
\begin{equation}
    \norm{\nabla f(\vect{x}_t+\vect{s}_t)}-2\epsilon_\vect{g} \geq \epsilon-2\epsilon_\vect{g}\geq \frac{1}{3}\epsilon
\end{equation}
Hence, from Eq. \ref{eq:s_t_bound_nabla} and Eq. \ref{eq:epsb_eps},
\begin{align}
    \norm{\vect{s}_t} &\geq \frac{1}{L_2+2\sigma_{\max}} \left( -\epsilon_\vect{B}+\sqrt{\epsilon_\vect{B}^2+\frac{2L_2+4\sigma_{\max}}{3}\epsilon} \right)\\
    & \geq \frac{1}{\sqrt{6}\sqrt{L_2+2\sigma_{\max}}}\epsilon^{\frac{1}{2}} \label{eq:s_epsilon_bound}
\end{align}
For a successful iterate $\vect{x}_t$, plugging Eq. \ref{eq:s_epsilon_bound} into Eq. \ref{eq:cubic_model_bound} and using Eq. \ref{eq:rho} yield,
\begin{align}
f(\vect{x}_t) - f(\vect{x}_{t+1}) &\geq \frac{\eta_1\sigma_t}{6}\norm{\vect{s}_t}^3   \\
& \geq \frac{\eta_1\sigma_{\min}}{6}\norm{\vect{s}_t}^3 \\
& \geq \frac{\eta_1\sigma_{\min}}{36\sqrt{6}\left( L_2+2\sigma_{\max} \right)^{\frac{3}{2}}}\epsilon^{\frac{3}{2}} = \kappa_{s1}\epsilon^{\frac{3}{2}}
\end{align}

From Lemma \ref{lemma:sufficient_decrease_nc} with $\epsilon'=\epsilon^{\frac{1}{2}}$ for an iteration $t \in \mathcal{U}_j$,
\begin{align}
f(\vect{x}_0)-\mathbb{E}[f(\vect{x}_j)] &= \sum_{t=0,t\in\mathcal{N}_j}^{j-1}\left[ f(\vect{x}_t)-f(\vect{x}_{t+1}) \right]+\sum_{t=0,t\in\mathcal{U}_j}^{j-1}\left[ f(\vect{x}_t)-\mathbb{E}[f(\vect{x}_{t+1})] \right] \label{eq:summing_up1}\\
&\geq |\mathcal{N}_j|\kappa_{s1}\epsilon^{\frac{3}{2}} + |\mathcal{U}_j|\kappa_{d}\epsilon^{\frac{3}{2}}\\
&\geq \left(|\mathcal{N}_j|+|\mathcal{U}_j|\right)\min\{\kappa_{s1}\epsilon^{\frac{3}{2}},\kappa_d\epsilon^{\frac{3}{2}}\}
\end{align}
That means 
\begin{equation}
    |\mathcal{N}_j|+|\mathcal{U}_j| \leq \frac{f(\vect{x}_0)-f_{low}}{\min\{\kappa_{s1}\epsilon^{\frac{3}{2}},\kappa_d\epsilon^{\frac{3}{2}}\}}
\end{equation}
with high probability 1-$\delta$.

Therefore, we can attain an upper bound on $j$,
\begin{equation}\label{eq:j_bound_1}
    j=|\mathcal{N}_j|+|\mathcal{U}_j| \leq \ceil[\bigg]{\frac{f(\vect{x}_0)-f_{low}}{\min\{\kappa_{s1}\epsilon^{\frac{3}{2}},\kappa_d\epsilon^{\frac{3}{2}}\}}}=l_1
\end{equation}

\item \label{enum:secondorder}
We also have to attain an upper bound for the second-order optimality as the same approach in \ref{enum:firstorder}.
For a sufficiently large iteration $j$, suppose that $-\lambda_{\min}(\vect{Q}_t^T\vect{B}_{t}\vect{Q}_t) > \epsilon , \,\forall t\leq j$.
Eq. \ref{eq:s_condition2} implies that the matrix $\vect{Q}_t^T\vect{B}_t\vect{Q}_t+\sigma_t\norm{\vect{s}_t}\vect{I}$ is positive semidefinite which gives
\begin{equation}
    \sigma_t\norm{\vect{s}_t} \geq |\lambda_{\min}(\vect{Q}_t^T\vect{B}_{t}\vect{Q}_t)| \text{ for all }t\leq j
\end{equation}
because $-\lambda_{\min}(\vect{Q}_t^T\vect{B}_{t}\vect{Q}_t) > \epsilon , \,\forall t\leq j$.\\
Hence, given an iteration $t\leq j$, from Eq. \ref{eq:cubic_model_bound} of Lemma \ref{lemma:local_cubic_model_decrease} and Eq. \ref{eq:rho} of $\rho_t$
\begin{align}
    f(\vect{x}_t) - f(\vect{x}_{t+1}) &\geq \frac{\eta_1\sigma_t}{6}\norm{\vect{s}_t}^3   \\
    & \geq \frac{\eta_1\sigma_{\min}}{6}\norm{\vect{s}_t}^3 \\
    & \geq \frac{\eta_1\sigma_{\min}}{6} \left(\frac{|\lambda_{\min}(\vect{Q}_t^T\vect{B}_{t}\vect{Q}_t)|}{\sigma_t}\right)^3\\
    & \geq \frac{\eta_1\sigma_{\min}}{6\sigma_{\max}}\epsilon^3 = \kappa_{s2}\epsilon^3
\end{align}
Summing up to iteration $j$ gives,
\begin{align}
    f(\vect{x}_0)-\mathbb{E}[f(\vect{x}_j)] &= \sum_{t=0,t\in\mathcal{N}_j}^{j-1}\left[ f(\vect{x}_t)-f(\vect{x}_{t+1}) \right]+\sum_{t=0,t\in\mathcal{U}_j}^{j-1}\left[ f(\vect{x}_t)-\mathbb{E}[f(\vect{x}_{t+1})] \right] \\
    &\geq |\mathcal{N}_j|\kappa_{s2}\epsilon^3+|\mathcal{U}_j|\kappa_d\epsilon^{\frac{3}{2}}\\
    &\geq \left(|\mathcal{N}_j|+|\mathcal{U}_j|\right)\min\{\kappa_{s2}\epsilon^3,\kappa_d\epsilon^{\frac{3}{2}}\}
\end{align}
Thus, the upper bound on $j$ is
\begin{equation}\label{eq:j_bound_2}
    j\leq l_2
\end{equation}

\item \label{enum:sumup}
To get  an upper bound on all the possible iterations that occur either with $\norm{\nabla f(\vect{x}_{t^*})} > \epsilon$ or $-\lambda_{\min}(\vect{Q}_{t^*}^T\vect{B}_{t^*}\vect{Q}_{t^*}) > \epsilon$ , the sum of the bounds Eq. \ref{eq:j_bound_1} and Eq. \ref{eq:j_bound_2} gives,
\begin{equation}
    j\leq l_1+l_2
\end{equation}
That equivalently means that if $j$ exceeds $l_1+l_2+1$, there exists at least an iterations $t$ which necessarily satisfy the condition $\norm{\nabla f(\vect{x}_{t^*})} \leq \epsilon$ and $-\lambda_{\min}(\vect{Q}_{t^*}^T\vect{B}_{t^*}\vect{Q}_{t^*}) \leq \epsilon$.
\end{enumerate}
\end{proof}

\newpage
\textbf{More numerical results}

\begin{figure}[h!]
\centering
\captionsetup{justification=centering}
  \subcaptionbox{w8a dataset\\($n=49749$, $d=300$)}{\includegraphics[height=2.in]{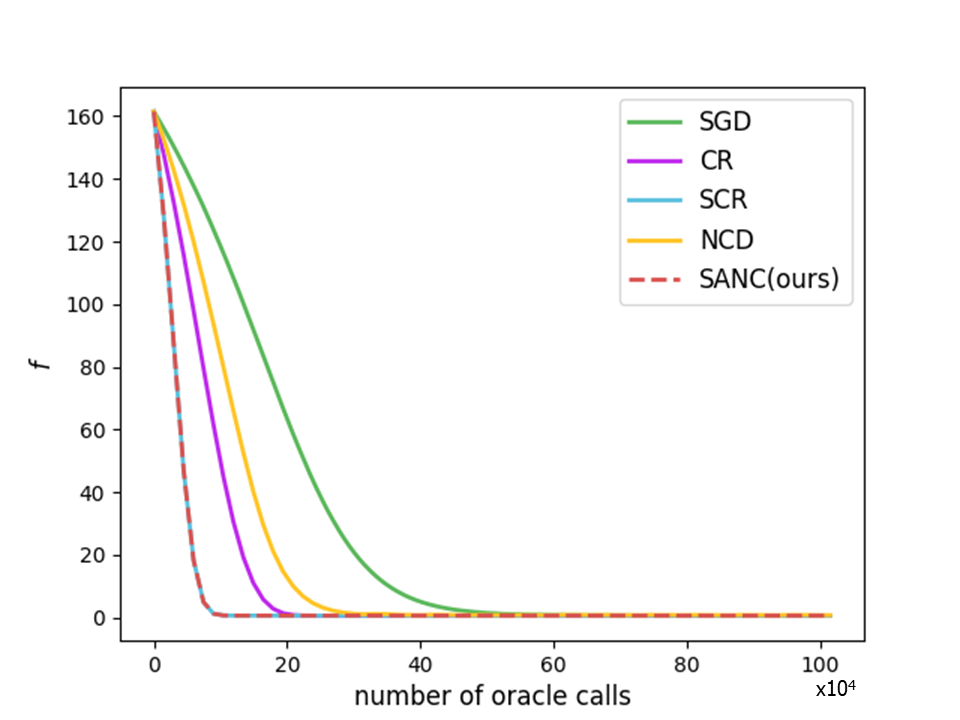}}\quad
  \subcaptionbox{a9a dataset\\($n=32561$, $d=123$)}{\includegraphics[height=2.in]{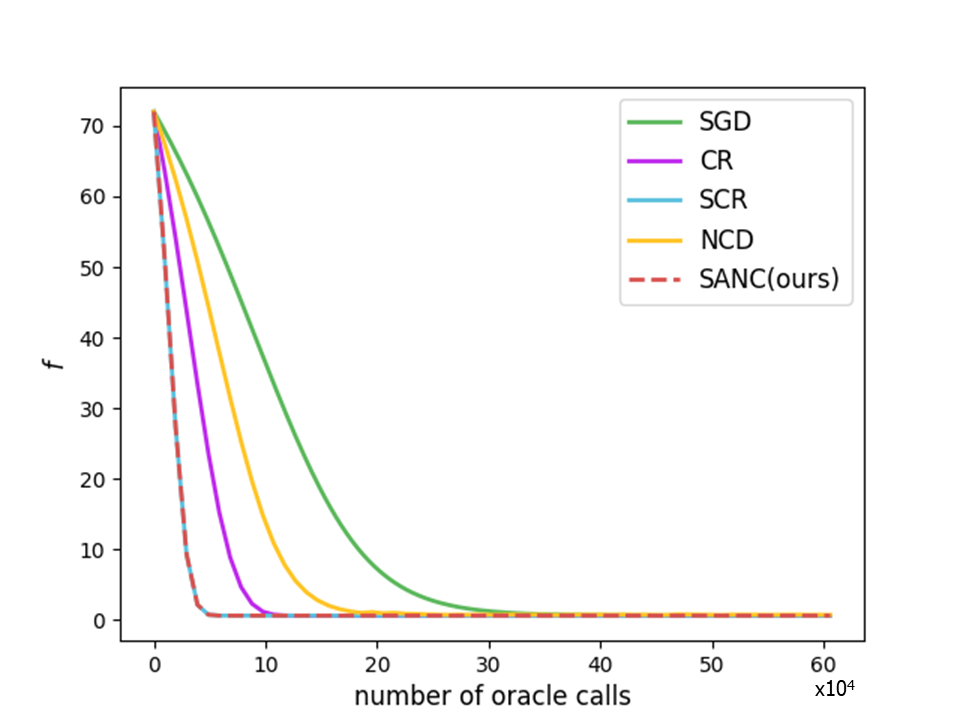}}\\
  \subcaptionbox{covtype dataset\\($n=581012$, $d=54$)}{\includegraphics[height=2.in]{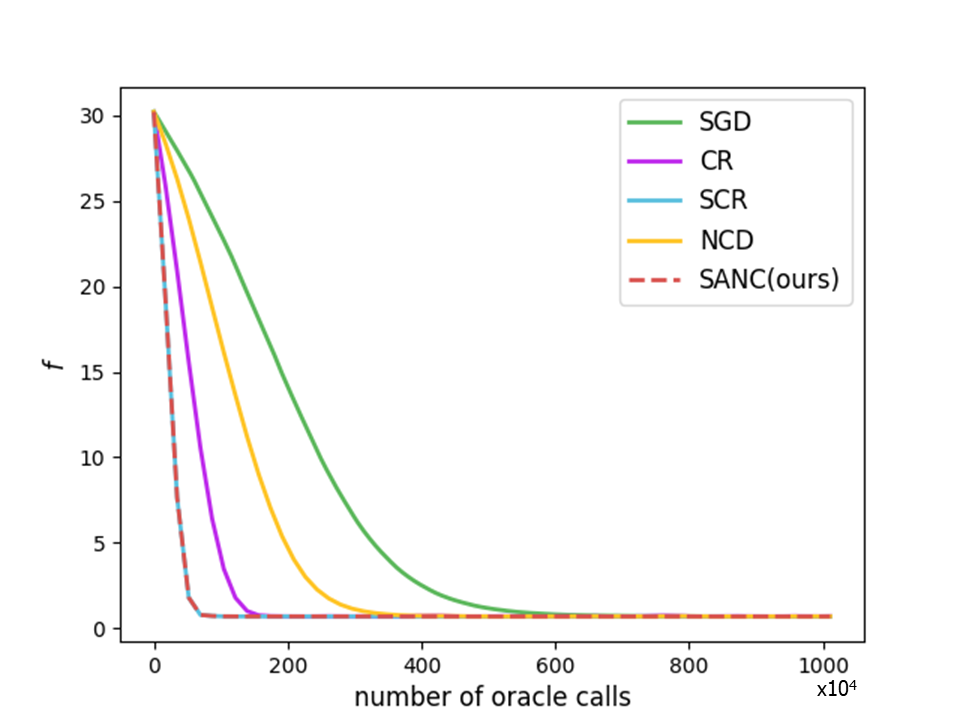}}\quad
  \subcaptionbox{ijcnn1 dataset\\($n=49990$, $d=22$)}{\includegraphics[height=2.in]{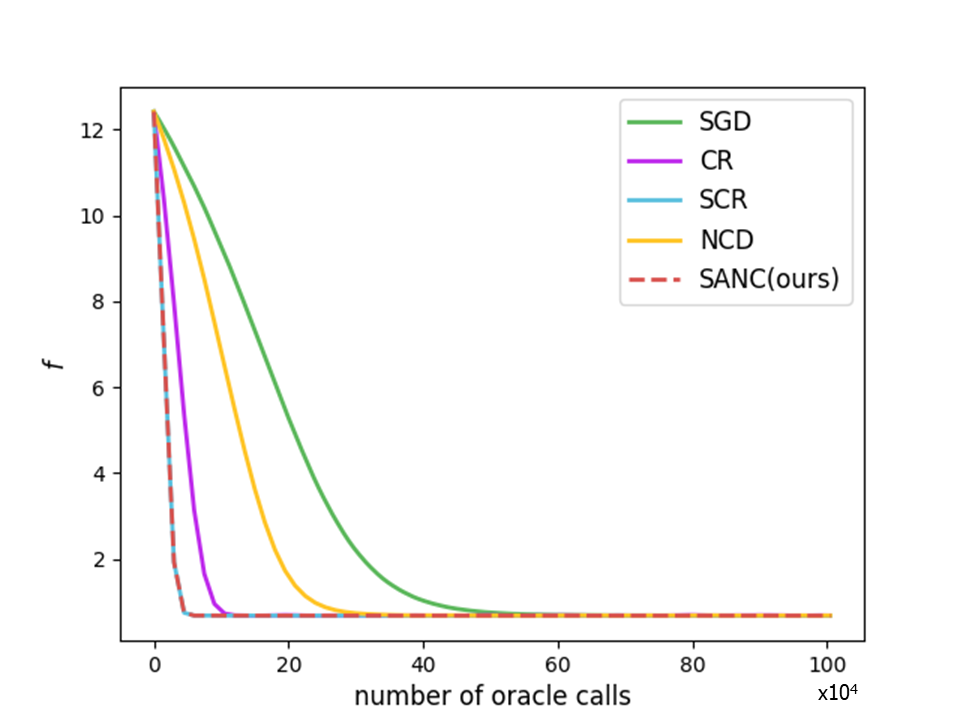}}
  \caption{Training losses of the logistic regression over the number of oracle calls. All function losses are the average of independent 10 runs.}\label{fig:lr_normal_more}
\end{figure}

% TODO: Figure 3-5를 appendix로 넣은 이유가?? 내가 보기엔 그냥 5. Empirical Studies로 옮기는게 낫겠는데 
\begin{figure*}[h!]
\begin{subfigure}{0.5\textwidth}
\centering
\captionsetup{justification=centering}
\includegraphics[height=2.7in]{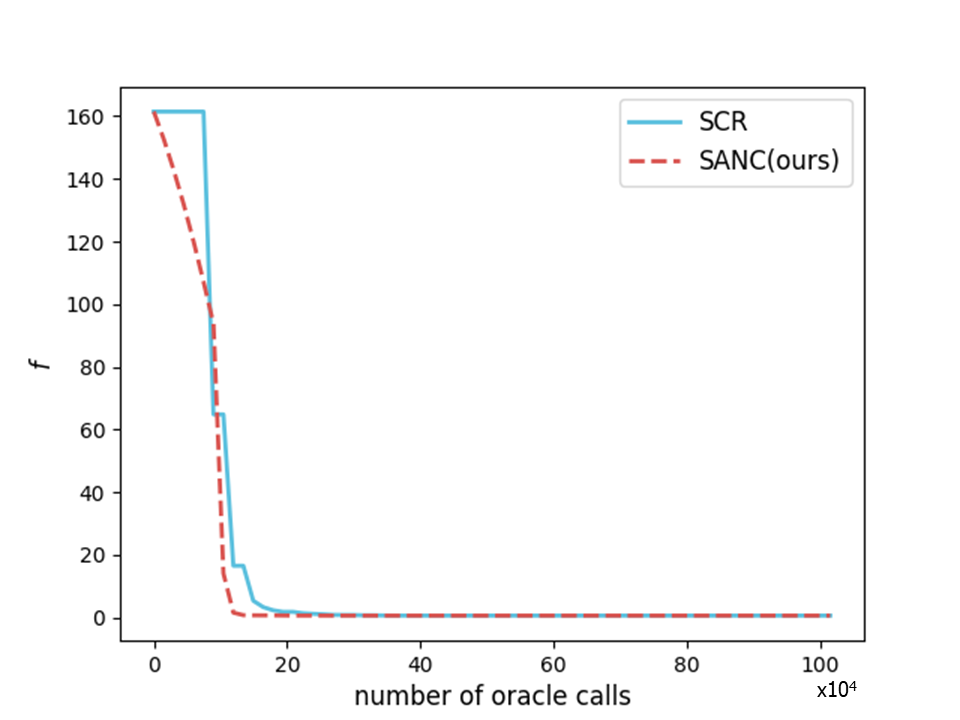}
\caption{w8a dataset ($n=49749$, $d=300$)}
\end{subfigure}\hfill
\begin{subfigure}{0.5\textwidth}
\centering
\captionsetup{justification=centering}
\includegraphics[height=2.7in]{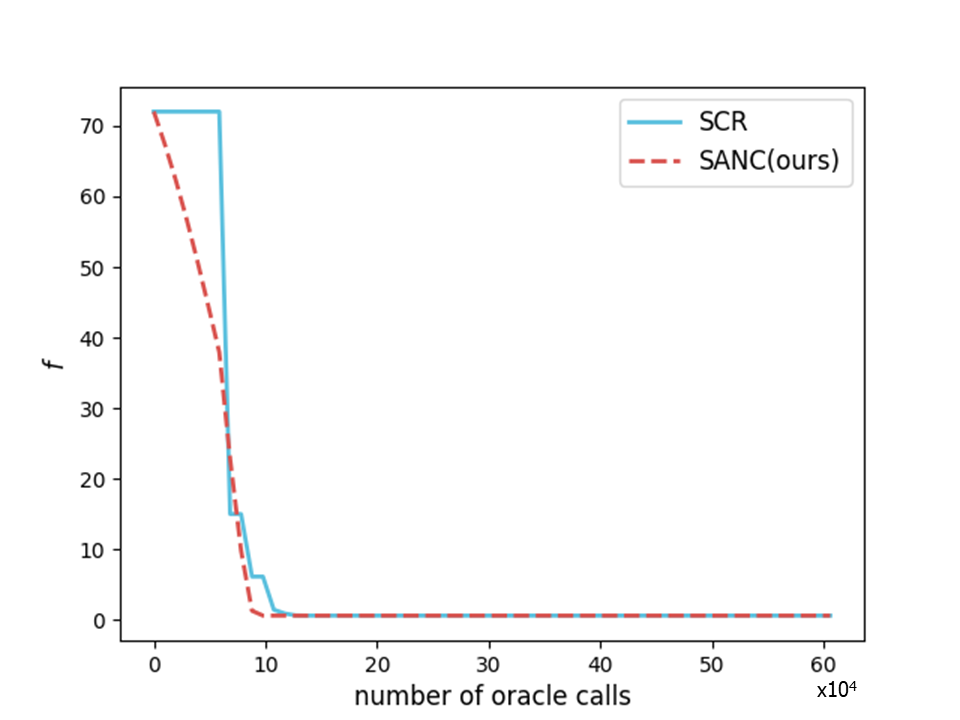}
\caption{a9a dataset ($n=32561$, $d=123$)}
\end{subfigure}
\caption{Training losses of the logistic regression with a nonconvex regularization term with $\lambda=1.0$ over the number of oracle calls. SANC and SCR start with $\sigma_0=0.001$. All function losses are the average of independent 10 runs.}
\end{figure*}

\begin{figure*}[h!]
\begin{subfigure}{0.5\textwidth}
\centering
\captionsetup{justification=centering}
\includegraphics[height=2.7in]{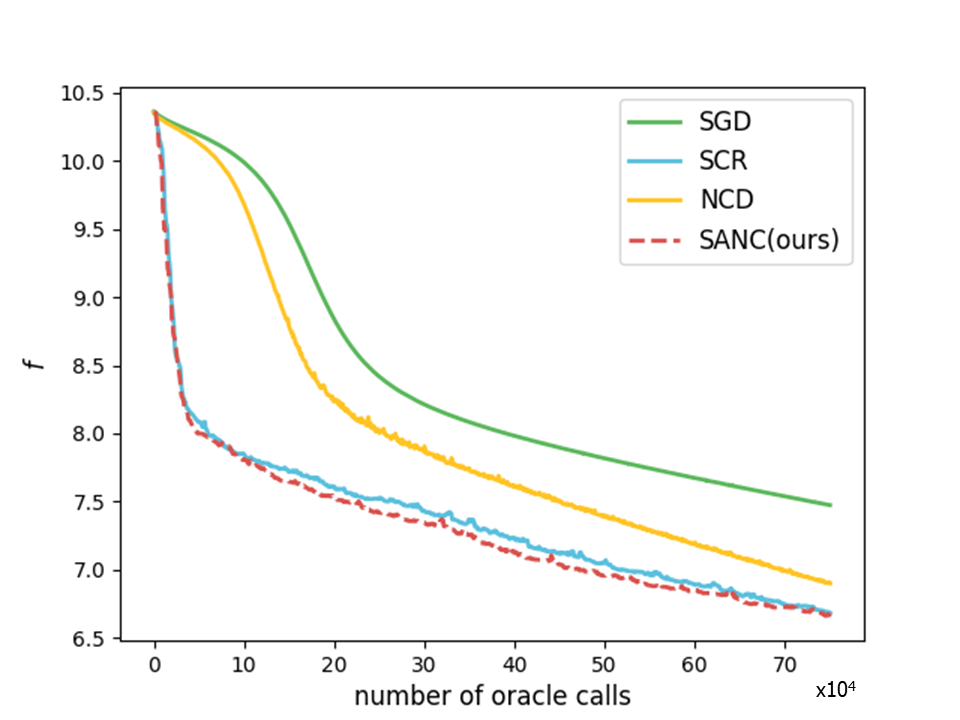}
\caption{MNIST dataset}
\end{subfigure}\hfill
\begin{subfigure}{0.5\textwidth}
\centering
\captionsetup{justification=centering}
\includegraphics[height=2.7in]{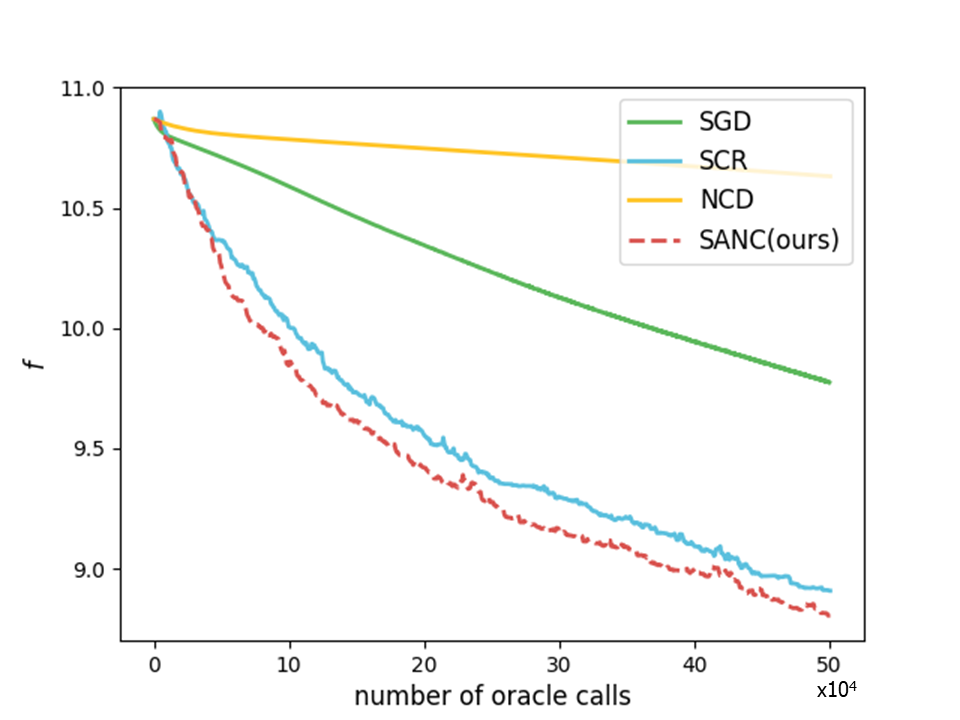}
\caption{CIFAR10 dataset}
\end{subfigure}
\caption{Training losses of the convolutional neural networks with ReLU nonlinear activation functions over the number of oracle calls. All function losses are the average of independent 10 runs.}\label{fig:CNN_relu}
\end{figure*}

\end{document}